\documentclass[12pt]{amsart}
\usepackage{xcolor}
\usepackage{mathtools}
\usepackage[latin2]{inputenc}
\usepackage{amsmath, amssymb}
\usepackage[
  hmarginratio={1:1},     
  vmarginratio={1:1},     
  textwidth=430pt,        
  heightrounded,          
]{geometry}

\usepackage{hyperref}

\newtheorem{theorem}{Theorem}
\newtheorem{lemma}{Lemma}
\newtheorem{proposition}{Proposition}
\newtheorem{algorithm}{Algorithm}

\theoremstyle{definition}
\newtheorem{definition}{Definition}

\theoremstyle{remark}
\newtheorem{remark}{Remark}

\newcommand{\esssup}{\operatorname{esssup}}
\newcommand{\Hm}{\mathcal{H}^{d-1}}
\newcommand{\bv}{BV}

\allowdisplaybreaks

\title[De Giorgi's inequality for the thresholding scheme]{De Giorgi's inequality for the thresholding scheme with arbitrary mobilities and surface tensions}

\author{Tim Laux}
\address[Tim Laux]{Institut f\"ur angewandte Mathematik, Universit\"at Bonn, Endenicher Allee 60, 53115 Bonn, Germany}
\email{tim.laux@hcm.uni-bonn.de}

\author{Jona Lelmi}
\address[Jona Lelmi]{Institut f\"ur angewandte Mathematik, Universit\"at Bonn, Endenicher Allee 60, 53115 Bonn, Germany}
\email{lelmi@hcm.uni-bonn.de}

\begin{document}

\maketitle

\begin{abstract}
	We provide a new convergence proof of the celebrated Merriman-Bence-Osher scheme for multiphase mean curvature flow. Our proof applies to the new variant incorporating a general class of surface tensions and mobilities, including typical choices for modeling grain growth. The basis of the proof are the minimizing movements interpretation of Esedo\u{g}lu and Otto and De Giorgi's general theory of gradient flows. Under a typical energy convergence assumption we show that the limit satisfies a sharp energy-dissipation relation.
\end{abstract}

\section{Introduction}

	The thresholding scheme is a highly efficient computational scheme for multiphase mean curvature flow (MCF) which was originally introduced by Merriman, Bence, and Osher \cite{Merriman1992}, \cite{Merriman1994}.
	The main motivation for MCF comes from metallurgy where it models the slow relaxation of grain boundaries in polycrystals \cite{mullins1956}. 
	Each ''phase'' in our mathematical jargon corresponds to a grain, i.e., a region of homogeneous crystallographic orientation. 
	The effective surface tension $\sigma_{ij}(\nu)$ and the mobility $\mu_{ij}(\nu)$ of a grain boundary depend on the mismatch between the lattices of the two adjacent grains $\Omega_i$ and $\Omega_j$ and on the relative orientation of the grain boundary, given by its normal vector $\nu$. It is well known that for small mismatch angles, the dependence on the normal can be neglected \cite{read1950dislocation}. 
	The effective evolution equations then read
	\begin{align}\label{eq:intro V=H}
		V_{ij} = -\mu_{ij} \sigma_{ij} H_{ij} \quad \text{along the grain boundary }\Sigma_{ij},
	\end{align}
	where $V_{ij} $ and $H_{ij}$ denote the normal velocity and mean curvature of the grain boundary $\Sigma_{ij} = \partial \Omega_i \cap \partial \Omega_j$, respectively. These equations are coupled by the Herring angle condition
	\begin{align}\label{eq:intro Herring}
		\sigma_{ij}\nu_{ij} + \sigma_{jk}\nu_{jk} + \sigma_{ki}\nu_{ki} =0 \quad \text{along triple junctions } \Sigma_{ij} \cap \Sigma_{jk},
	\end{align}
	which is a balance-of-forces condition and simply states that triple junctions are in local equilibrium. 
	Efficient numerical schemes allow to carry out large-scale simulations to give insight into relevant statistics like the average grain size or the grain boundary character distribution, as an alternative to studying corresponding mean field limits as in \cite{barmak2011critical}, \cite{klobusicky2020two}. 
	The main obstruction to directly discretize the dynamics \eqref{eq:intro V=H}--\eqref{eq:intro Herring} are ubiquitous topological changes in the network of grain boundaries like for example the vanishing of grains. 
	Thresholding instead naturally handles such topological changes. The scheme is a time discretization which alternates between the following two operations: \emph{(i)} convolution with a smooth kernel; \emph{(ii)} thresholding. The second step is a simple pointwise operation and also the first step can be implemented efficiently using the Fast Fourier Transform. 
	One of the main objectives of our analysis is to rigorously justify this intriguingly simple scheme in the presence of such topological changes.
	
	\medskip
	
	The basis of our analysis is the underlying gradient-flow structure of \eqref{eq:intro V=H}--\eqref{eq:intro Herring}, which means that the solution follows the steepest descent in an energy landscape. More precisely, the energy is the total interfacial area weighted by the surface tensions $\sigma_{ij}$, and the metric tensor is the $L^2$-product on normal velocities, weighted by the inverse mobilities $\frac1{\mu_{ij}}$.
	 One can read off this structure from the inequality
	 \begin{align*}
	 	\frac{d}{dt} \sum_{i,j=1}^N \sigma_{ij} \textup{Area}(\Sigma_{ij}) = - \sum_{i,j=1}^N \frac1{\mu_{ij}} \int_{\Sigma_{ij}} V_{ij}^2\, dS \leq 0,
	 \end{align*}
	 which is valid for sufficiently regular solutions to \eqref{eq:intro V=H}--\eqref{eq:intro Herring}. 
	In the seminal work \cite{Esedoglu2015}, Esedo\u{g}lu and Otto showed that the efficient thresholding scheme respects this gradient-flow structure as it may be viewed as a minimizing movements scheme in the sense of De Giorgi. More precisely, they show that each step in the scheme is equivalent to solving a variational problem of the form
	\begin{align}\label{eq:intro min mov}
		\min_\chi  \left\{ \frac1{2h} d_h^2(\Sigma,\Sigma^{n-1}) + E_h(\Sigma) \right\},
	\end{align}
	where $E_h(\Sigma)$ and $d_h(\Sigma,\Sigma^{n-1})$ are proxies for the total interfacial energy of the configuration $\Sigma$ and the distance of the configuration $\Sigma$ to the one at the previous time step $\Sigma^{n-1}$, respectively. Since the work of Jordan, Kinderlehrer, and Otto \cite{Jordan1998}, the importance of the formerly often neglected metric in such gradient-flow structures has been widely appreciated.  
	Also in the present work, the focus lies on the metric, which in the case of MCF is well-known to be completely degenerate \cite{Michor2006}. This explains the proxy for the metric appearing in the related well-known minimizing movements scheme for MCF by Almgren, Taylor, and Wang, \cite{Almgren1993}, and Luckhaus and Sturzenhecker \cite{Luckhaus1995}. 
	This remarkable connection between the numerical scheme and the theory of gradient flows has the practical implication that it made clear how to generalize the algorithm to arbitrary surface tensions $\sigma_{ij}$. 
	From the point of view of numerical analysis, \eqref{eq:intro min mov} means that thresholding behaves like the \emph{implicit} Euler scheme and is therefore unconditionally stability. 
	The variational interpretation of the thresholding scheme has of course implications for the analysis of the algorithm as well. It allowed Otto and one of the authors to prove convergence results in the multiphase setting \cite{Laux2016}, \cite{Laux2020}, which lies beyond the reach of the more classical viscosity approach based on the comparison principle implemented in \cite{Evans1993}, \cite{Barles1995}, \cite{Ishii1999}. Also in different frameworks, this variational viewpoint turned out to be useful, such as MCF in higher codimension \cite{Laux2019a} or the Muskat problem \cite{Jacobs2021}. 
	The only downside of the generalization \cite{Esedoglu2015} are the somewhat unnatural effective mobilities $\mu_{ij}=\frac1{\sigma_{ij}}$. 
	Only recently, Salvador and Esedo\u{g}lu \cite{Salvador2019} have presented a strikingly simple way to incorporate a wide class of mobilities $\mu_{ij}$ as well. 
	Their algorithm is based on the fact that although the same kernel appears in the energy and the metric, each term only uses certain properties of the kernel, which can be tuned independently: Starting from two Gaussian kernels $G_\gamma$ and $G_\beta$ of different width, they find a positive linear combination $K_{ij}=a_{ij} G_\gamma + b_{ij} G_\beta$, whose effective mobility and surface tension match the given $\mu_{ij}$ and $\sigma_{ij}$, respectively. It is remarkable that this algorithm retains the same simplicity and structure as the previous ones \cite{Merriman1994}, \cite{Esedoglu2015}. We refer to Section \ref{setup} for the precise statement of the algorithm.
	
	In the present work, we prove the first convergence result for this new general scheme. We exploit the gradient-flow structure and show that under the natural assumption of energy convergence, any limit of thresholding satisfies De Giorgi's inequality, a weak notion of multiphase mean curvature flow. This assumption is inspired by the fundamental work of Luckhaus-Sturzenhecker \cite{Luckhaus1995} and has appeared in the context of thresholding in \cite{Laux2016}, \cite{Laux2020}. 
	We expected it to hold true before the onset of singularities such as the vanishing of grains. Furthermore, at least in the simpler two-phase case, it can be verified for certain singularities \cite{Philippis2018}, \cite{Chambolle2006}. We would in fact expect this assumption to be true generically, which however seems to be a difficult problem in the multiphase case.
	
	
	The present work fits into the theory of general gradient flows even better than the two previous ones \cite{Laux2016}, \cite{Laux2020} and crucially depends on De Giorgi's abstract framework, cf.\ \cite{Ambrosio2005}. This research direction was initiated by Otto and the first author and appeared in the lecture notes \cite{Laux2019}. There, De Giorgi's inequality is derived for the simple model case of two phases. Here, we complete these ideas and use a careful localization argument to generalize this result to the multiphase case. A further particular novelty of our work is that for the first time, we prove the convergence of the new scheme for arbitrary mobilities \cite{Salvador2019}.
	
	\medskip
	
	Our proof rests on the fact that thresholding, like any minimizing movements scheme, satisfies a sharp energy-dissipation inequality of the form
	\begin{align}\label{eq:intro EDI}
		E_h(\Sigma^h(T)) +  \frac12 \int_0^T \left( \frac1{h^2}d_h^2(\Sigma^h(t),\Sigma^h(t-h))  
		+|\partial E_h|^2(\tilde \Sigma^h(t))  \right) dt \leq E_h(\Sigma(0)),
	\end{align}
	where $\Sigma^h(t)$ denotes the piecewise constant interpolation in time of our approximation,   $\tilde \Sigma^h(t)$ denotes another, intrinsic interpolation in terms of the variational scheme, cf.\ Lemma \ref{afpLemma}, and $|\partial E_h|$ is the metric slope of $E_h$, cf.\ (\ref{metricDerDef}). 
	
	Our main goal is to pass to the limit in \eqref{eq:intro EDI} and obtain the sharp energy-dissipation relation for the limit, which in the simple two-phase case formally reads
	\begin{align}\label{eq:intro DeGiorgi}
		\sigma \mathrm{Area}(\Sigma(T)) + \frac1{2} \int_0^T \int_{\Sigma(t)} \left( \frac1\mu V^2
		+ \sigma^2 \mu H^2 \right) dS\,dt \leq \sigma \mathrm{Area}(\Sigma(0)).
	\end{align}	
	To this end, one needs sharp lower bounds for the terms on the left-hand side of \eqref{eq:intro EDI}.  
	While the proof of the lower bound on the metric slope of the energy
	\begin{align}\label{eq:intro LSC H}
		 \liminf_{h\downarrow0} \int_0^T  |\partial E_h|^2(\tilde\Sigma^h(t)) \,dt 
		 \geq 	\sigma^2 \mu \int_0^T \int_{\Sigma(t)} H^2 dS\,dt
	\end{align}
	is a straight-forward generalization of the argument in \cite{Laux2019}, the main novelty of the present work lies in the sharp lower bound for the distance-term of the form
	\begin{align}\label{eq:intro LSC V}
		\liminf_{h\downarrow0} \int_0^T \frac1{h^2} d_h^2(\Sigma^h(t),\Sigma^h(t-h))\,dt
		\geq \frac1\mu \int_0^T\int_{\Sigma(t)} V^2\, dS\,dt.
	\end{align}
	This requires us to work on a mesoscopic time scale $\tau \sim \sqrt{h}$, which is much larger than the microscopic time-step size $h$ and which is natural in view of the parabolic nature of our problem. 
	It is remarkable that De Giorgi's inequality \eqref{eq:intro DeGiorgi} in fact characterizes the solution of MCF under additional regularity assumptions. 
	Indeed, if $\Sigma(t)$ evolves smoothly, this inequality can be rewritten as
	\begin{align}
	\frac1{2} \int_0^T \int_{\Sigma(t)}
	\sigma \Big( \frac1{\sqrt{\mu\sigma}} V
	+ \sqrt{\sigma \mu} H \Big)^2 dS\,dt \leq 0,
	\end{align}
	and therefore $V=-\mu \sigma H$. 
	For expository purpose, we focused here on the vanilla two-phase case. In the multiphase case, the resulting inequality implies both the PDEs \eqref{eq:intro V=H} and the balance-of-forces conditions \eqref{eq:intro Herring}, cf.\ Remark \ref{relationSmoothSol}.
	An optimal energy-dissipation relation like the one here also plays a crucial role in the recent weak-strong uniqueness result for multiphase mean curvature flow by Fischer, Hensel, Simon, and one of the authors \cite{Fischer2020}. There, a new dynamic analogue of calibrations is introduced and uniqueness is established in the following two steps:  \emph{(i)} any strong solution is a calibrated flow and \emph{(ii)} every calibrated flow is unique in the class of weak solutions.
	De Giorgi's general strategy we are implementing here is also related to the approaches by Sandier and Serfaty \cite{Sandier2004} and Mielke \cite{mielke2016}. They provide sufficient conditions for gradient flows to converge in the same spirit as $\Gamma$-convergence of energy functionals, implies the convergence of minimizers. In the dynamic situation it is clear that one needs conditions on both energy and metric in order to verify such a convergence.
	
	\medskip
	
	There has been continuous interest in MCF in the mathematics literature, so we only point out some of the most relevant recent advances. We refer the interested reader to the introductions of \cite{Laux2016} and \cite{Laux2018} for further related references. 
	The existence of global solutions to multiphase MCF has only been established recently by Kim and Tonegawa \cite{Kim2017} who carefully adapt Brakke's original construction and show in addition that phases do not vanish spontaneously. For the reader who wants to familiarize themselves with this topic, we recommend the recent notes \cite{Tonegawa2019}. 
	Another approach to understanding the long-time behavior of MCF flow is to restart strong solutions after singular times. This amounts to solving the Cauchy problem with non-regular initial data, such as planar networks of curves with quadruple junctions. 
	In this two-dimensional setting, this has been achieved by Ilmanen, Neves, and Schulze \cite{Ilmanen2019} by gluing in self-similarly expanding solutions for which it is possible to show that the initial condition is attained in some measure theoretic way. 
	Most recently, using a similar approach of gluing in self-similar solutions, but also relying on blow-ups from geometric microlocal analysis, Lira, Mazzeo, Pluda, Saez \cite{Lira2021} were able to construct such strong solutions, prove stronger convergence towards the initial (irregular) network of curves, and classify all such strong solutions.
	
	\medskip
	
	The rest of the paper is structured as follows. In Section \ref{setup} we recall the thresholding scheme for arbitrary mobilitites introduced in \cite{Salvador2019}, show its connection to the abstract framework of gradient flows, and record the direct implications of this theory. We state and discuss our main results in Section \ref{sec:results}.  Section \ref{consPOU} contains the localization argument in space, which will play a crucial role in the proofs which are gathered in Section \ref{sec:proofs}. Finally, in the short Appendix, we record some basic facts about thresholding.
\section{Setup and the modified thresholding scheme}\label{setup}

\subsection{The modified algorithm}

We start by describing the algorithm proposed by Salvador and Esedo\u{g}lu in \cite{Salvador2019}. Let the symmetric matrix $\mathbb{\sigma} = (\sigma_{ij})_{ij} \in \mathbf{R}^{N\times N}$ of surface tensions and the symmetric matrix $\mathbb{\mu} = (\mu_{ij})_{ij}$ of mobilities be given. In this work we define for notational convenience $\sigma_{ii} = \mu_{ii} = 0$. Let $\gamma > \beta > 0$ be given. Define the matrices $\mathbb{A} = (-a_{ij})_{ij} \in \mathbf{R}^{N \times N}$ and $\mathbb{B} = (-b_{ij})_{ij} \in \mathbf{R}^{N \times N}$ by

\begin{align}
& a_{ij} = \frac{\sqrt{\pi}\sqrt{\gamma}}{\gamma - \beta}(\sigma_{ij} - \beta \mu_{ij}^{-1}),
\\ & b_{ij} = \frac{\sqrt{\pi}\sqrt{\beta}}{\gamma - \beta}(-\sigma_{ij} + \gamma \mu_{ij}^{-1}),
\end{align}
for $i \neq j$ and $a_{ii} = b_{ii} = 0$. Then $a_{ij}, b_{ij}$ are uniquely determined as solutions of the following linear system

\begin{equation}
\begin{cases}
\sigma_{ij} = \frac{a_{ij}\sqrt{\gamma}}{\sqrt{\pi}} + \frac{b_{ij}\sqrt{\beta}}{\sqrt{\pi}},
\\
\mu_{ij}^{-1} = \frac{a_{ij}}{\sqrt{\pi}\sqrt{\gamma}} + \frac{b_{ij}}{\sqrt{\pi}\sqrt{\beta}}.
\end{cases}
\end{equation}
The algorithm introduced by Salvador and Esedo\u{g}lu is as follows. Let the time step size $h > 0$ be fixed. Hereafter $G_{\gamma}^h := G_{\gamma h}^{(d)}$ denotes the $d$-dimensional heat kernel (\ref{heatK}) at time $\gamma h$.
\begin{algorithm}[Modified thresholding scheme]\label{algo}
Given the initial partition $\Omega_1^0, . . . , \Omega_N^0$, to obtain the partition $\Omega_1^{n+1}, . . . , \Omega_N^{n+1}$ at time $t= h (n+1)$ from the partition $\Omega_1^n, . . . , \Omega_N^n$ at time $t = hn$

\begin{enumerate}
\item For any $i = 1 , . . . , N$ form the convolutions

\begin{equation}
\phi^n_{1,i} = G_{\gamma}^h * \mathbf{1}_{\Omega_i^n},\ \phi^n_{2,i} = G_{\beta}^h * \mathbf{1}_{\Omega_i^n}\nonumber
\end{equation}

\item For any $i= 1, . . . , N$ form the comparison functions

\begin{equation}\label{secondStepT}
\psi^n_{i} = \sum_{j\neq i} a_{ij}\phi^n_{1,j} + b_{ij}\phi^n_{2,j}.\nonumber
\end{equation}

\item Thresholding step, define

\begin{equation}
\Omega_i^{n+1} := \left\{x: \psi_i^n(x) < \min_{j \neq i} \psi_j^n(x) \right\}.\nonumber
\end{equation}
\end{enumerate}
\end{algorithm}
We will assume the following:

\begin{align}
&\text{The coefficients}\ a_{ij}, b_{ij}\ \text{satisfy the strict triangle inequality.}\label{strictTI}
\\ 
&\text{The matrices}\ \mathbb{A}\ \text{and}\ \mathbb{B}\ \text{are positive definite on}\ (1, . . . , 1)^{\perp}.\label{posDefAss}
\end{align}
In particular, for $v \in (1, . . . , 1)^{\perp}$ we can define norms

\begin{align}\label{normDef}
|v|_{\mathbb{A}}^2 = v \cdot \mathbb{A}v,\ |v|_{\mathbb{B}}^2 = v \cdot \mathbb{B}v.\nonumber
\end{align}

Observe that condition (\ref{strictTI}) is always satisfied if we choose $\gamma$ large and $\beta$ small provided the surface tensions and the inverse of the mobilities satisfy the strict triangle inequality. Indeed define

\begin{align*}
m_{\sigma} = \min_{i,j,k} \{\sigma_{ik} + \sigma_{kj} - \sigma_{ij}\}\ \text{and}\ M_{\sigma} = \max_{i,j,k} \{\sigma_{ik} + \sigma_{kj} - \sigma_{ij}\},
\end{align*}
where $i,j,k$ range over all triples of distinct indices $1\le i,j,k,\le N$. Define $m_{\frac{1}{\mu}}$ and $M_{\frac{1}{\mu}}$ in a similar way. Then a computation shows that $a_{ij}$ and $b_{ij}$ satisfy the (strict) triangle inequality if

\begin{equation}
\beta < \frac{m_{\sigma}}{M_{\frac{1}{\mu}}}\ \text{and}\ \gamma > \frac{M_{\sigma}}{m_{\frac{1}{\mu}}},
\end{equation}
which can always be achieved since $\gamma > \beta > 0$ are arbitrary. For the second condition (\ref{posDefAss}), we have the following result of Salvador and Esedo\u{g}lu \cite{Salvador2019}.

\begin{lemma}\label{PosDef}
Let the matrix $\sigma$ of the surface tensions and the matrix $\frac{1}{\mu}$ of the inverse mobilities (for the diagonal we set inverses to be zeros) be negative definite on $(1, . . . , 1)^\perp$. Let $\gamma > \beta$ be such that

\begin{equation}
\gamma > \frac{\min_{i=1, . . . , N-1} s_i}{\max_{i=1, . . . , N-1} m_i},\ \beta < \frac{\max_{i=1, . . . , N-1}s_i}{\min_{i=1, . . . , N-1}m_i}
\end{equation}
where $s_i$ and $m_i$ are the nonzero eigenvalues of $J\sigma J$ and $J\frac{1}{\mu}J$ respectively, where the matrix $J$ has components $J_{ij} = \delta_{ij} - \frac{1}{N}$. Then $\mathbb{A}$ and $\mathbb{B}$ are positive definite on $(1, . . . , 1)^{\perp}$.
\end{lemma}
In particular, if we choose $\gamma$ large enough and $\beta$ small enough, condition (\ref{posDefAss}) on the matrices $\mathbb{A}, \mathbb{B}$ is satisfied provided the matrices $\sigma$ and $\frac{1}{\mu}$ are negative definite on $1, . . . , 1)^{\perp}$. By a classical result of Schoenberg \cite{Schoenberg1938} this is the case if and only if $\sqrt{\sigma_{ij}}$ and $1/\sqrt{\mu_{ij}}$ are $\ell^2$ embeddable. In particular, this holds for the choice of Read-Schockley surface tensions and equal mobilities.

For $1 \le i \neq j \le N$ define the kernels 

\begin{equation}\label{kernelDef}
K_{ij}(z) = a_{ij}G_{\gamma}(z) + b_{ij}G_{\beta}(z)
\end{equation}
where, for a given $t > 0$, we define $G^{(d)}_t$ as the heat kernel in $\mathbf{R}^d$, i.e.,

\begin{equation}\label{heatK}
G^{(d)}_t(z) = \frac{e^{-\frac{|z|^2}{4t}}}{\sqrt{4\pi t}^d}.
\end{equation}
If the dimension $d$ is clear from the context, we suppress the superscript $(d)$ in (\ref{heatK}). We recall here some basic properties of the heat kernel.

\begin{align}
&G_t(z) > 0\ \text{(non-negativity)},\label{nonnegHK}
\\ &G_t(z) = G_t(Rz)\ \forall R \in O(d)\ \text{(symmetry)},\label{symmetryHK}
\\ &G_{t}(z) = \frac{1}{\sqrt{t}^d}G_1\left(\frac{z}{\sqrt{t}}\right)\ \text{(scaling)},\label{scalingHK}
\\ &G_t * G_s = G_{t+s}\ \text{(semigroup property)},\label{semigroupHK}
\\ &G_t^{(d)}(z) = \prod_{i=1}^d G_t^{(1)}(z_i)\ \text{(factorization property)}.\label{factpropHK} 
\end{align}
We observe that the kernels $K_{ij}$ are positive, with positive Fourier transform $\hat{K}_{ij}$ provided $\gamma > \max_{i,j}\sigma_{i,j}\mu_{i,j}$ and $\beta < \min_{i,j}\sigma_{i,j}\mu_{i,j}$. In particular assuming

\begin{enumerate}
\item $\sigma_{ij}$ and $\frac{1}{\mu_{ij}}$ satisfy the strict triangle inequality,
\item $\sigma$ and $\frac{1}{\mu}$ are negative definite on $(1, . . . , 1)^{\perp}$,
\end{enumerate}
we can always achieve the conditions posed on $\mathbb{A}, \mathbb{B}$ and the positivity of the kernels $K_{ij}$ by choosing $\gamma$ large and $\beta$ small.

Given any $h > 0$ we define the scaled kernels

\begin{equation}\label{scaledK}
K^h_{ij}(z) = \frac{1}{\sqrt{h}^{d}}K_{ij}(\frac{z}{\sqrt{h}}),
\end{equation}
then the first and the second step in Algorithm \ref{algo} may be compactly rewritten as follows
\begin{align*}
\psi_i^n = \sum_{j\neq i} K_{ij}^h * \mathbf{1}_{\Omega_j^n}.
\end{align*}
For later use, we also introduce the kernel

\begin{equation}
K(z) = \frac{1}{2}G_{\gamma}(z) + \frac{1}{2}G_{\beta}(z).
\end{equation}

\subsection{Connection to De Giorgi's minimizing movements}

The  first observation is that Algorithm \ref{algo} has a minimizing movements interpretation. To explain this, let us introduce the class

\begin{equation}
\mathcal{A} := \left\{ \chi: [0,1)^d \to \{0,1\}^N\biggr\rvert\ \sum_{k=1}^N \chi_k = 1 \right\}\nonumber
\end{equation}
and its relaxation

\begin{equation}
\mathcal{M} := \left\{ u: [0,1)^d \to [0,1]^N\biggr\rvert\ \sum_{k=1}^N u_k = 1 \right\}.\nonumber
\end{equation}
If $\chi \in \mathcal{A}\cap BV([0,1)^d)^N$, then each of the sets $\Omega_i := \{\chi_i = 1\}$ is a set of finite perimeter. We denote by $\partial^*\Omega_i$ the reduced boundary of the set $\Omega_i$, and for any pair $1\le i\neq j\le N$ we denote by $\Sigma_{ij} := \partial^*\Omega_i \cap \partial^*\Omega_j$ the interface between the sets. For $u \in \mathcal{M}$ we define

\begin{equation}\label{energyDef}
E(u) := \begin{cases}
\sum_{i,j}\sigma_{ij}\Hm(\Sigma_{ij})\ &\text{if}\ u \in \mathcal{A} \cap BV([0,1)^d)^N
\\
+\infty\ &\text{otherwise}.
\end{cases}
\end{equation}
For $h>0$ fixed we define the approximate energy $E_h$ for $u \in \mathcal{M}$

\begin{equation}\label{energyApprox}
E_h(u) = \sum_{i,j} \frac{1}{\sqrt{h}}\int_{[0,1)^d} u_iK_{ij}^h * u_j dx.
\end{equation}
For $u,v \in \mathcal{M}$ and $h>0$ we also define the distance

\begin{align}\label{distanceDef}
d_h^2(u,v) &:= -2hE_h(u-v) = -2\sqrt{h}\sum_{i,j} \int (u_i - v_i)K_{ij}^h*(u_j - v_j)dx
\\ &=2\sqrt{h}\int |G_{\gamma}^{h/2} *(u-v)|_{\mathbb{A}}^2 + |G_{\beta}^{h/2} *(u-v)|_{\mathbb{B}}^2\ dx,\nonumber
\end{align}
where we used the semigroup property (\ref{semigroupHK}) and the symmetry (\ref{symmetryHK}) to derive the last equality.

\begin{lemma}\label{minmovInt}
The pair $(\mathcal{M}, d_h)$ is a compact metric space. The function $E_h$ is continuous with respect to $d_h$. For every $1 \le i \le N$ and $n \in \mathbf{N}$ define $\chi_i^n = \mathbf{1}_{\Omega^n_i}$, where $\Omega_1^n, . . . , \Omega_N^n$ are obtained from $\Omega_1^{n-1}, . . . , \Omega_N^{n-1}$ by the thresholding scheme. Then $\chi^n$ minimizes
\begin{equation}\label{whatMinimizes}
\frac{1}{2h}d_h^2(u, \chi^{n-1}) + E_h(u)\ \text{among all}\ u \in \mathcal{M}.
\end{equation}
\end{lemma}

\begin{proof}
For $u, v \in \mathcal{M}$ definition (\ref{distanceDef}) and the fact that $\mathbb{A}$ and $\mathbb{B}$ are positive definite imply that $d_h$ is a distance on $\mathcal{M}$. The fact that $(\mathcal{M}, d_h)$ is compact and $E_h$ is continuous is just a consequence of the fact that $d_h$ metrizes the weak convergence in $L^2$ on $\mathcal{M}$, the interested reader may find the details of the reasoning in \cite{Laux2019}. We are thus left with showing that $\chi^n$ satisfies (\ref{whatMinimizes}). For $u, v \in L^2([0,1)^d)$ define

\begin{equation}
(u,v) = \frac{1}{\sqrt{h}}\sum_{i,j}\int u_i K^{h}_{ij}*v_j dx,\nonumber
\end{equation}
then by the symmetry (\ref{symmetryHK}) of the Gaussian kernel and by the symmetry of both matrices $\mathbb{A}, \mathbb{B}$ it is not hard to show that $(\cdot, \cdot)$ is symmetric. In particular we can write for any $u \in \mathcal{M}$

\begin{align*}
\frac{1}{2h}d_h^2(u,\chi^{n-1}) + E_h(u) &= -E_h(u-\chi^{n-1}) + E_h(u)
\\ & =-(u-\chi^{n-1}, u-\chi^{n-1}) + (u,u)
\\ & =2(\chi^{n-1}, u) - (\chi^{n-1}, \chi^{n-1}).
\end{align*}
Thus (\ref{whatMinimizes}) is equivalent to the fact that $\chi^n$ minimizes $(\chi^{n-1}, u)$ among all $u \in \mathcal{M}$. Since by (\ref{secondStepT})

\begin{equation}
(\chi^{n-1}, u) =  \int \sum_{i} u_i \psi^n_i dx,\nonumber
\end{equation}
we see that $\chi^{n}$ minimizes the integrand pointwise, and thus it is a minimizer for the functional.\end{proof}
The previous lemma allows us to apply the general theory of gradient flows in \cite{Ambrosio2005} to this particular problem. We record the key statement for our purposes in the following lemma.

\begin{lemma}\label{afpLemma}

Let $(\mathcal{M},d)$ be a compact metric space and $E:\mathcal{M} \to \mathbf{R}$ be continuous. Given $\chi^0 \in \mathcal{M}$ and $h > 0$ consider a sequence $\{\chi^n\}_{n\in \mathbf{N}}$ satisfying

\begin{equation}\label{whatchin}
\chi^n\ \text{minimizes}\ \frac{1}{2h}d^2(u, \chi^{n-1})+ E(u)\ \text{among all}\ u\in \mathcal{M}.
\end{equation}

Then we have for all $t \in \mathbf{N}h$

\begin{align}\label{DGineq}
\begin{aligned}
E(\chi(t)) + \frac{1}{2}\int_0^t\left(\frac{1}{h^2}d^2(\chi(s+h),\chi(s)) +|\partial E|^2(u(s))\right) ds \le E(\chi^0).
\end{aligned}
\end{align}

Here $\chi(t)$ is the piecewise constant interpolation, $u(t)$ is another interpolation satisfying

\begin{align}
&\begin{aligned}\label{varDist}
\int_0^{\infty}\frac{1}{2h^2}d^2(u(t),\chi(t))dt \le E(\chi^0),
\end{aligned}
\\ &\begin{aligned}\label{varEnergy}
\\ &E(u(t)) \le E(\chi(t))\ \text{for all}\ t\ge 0,
\end{aligned}
\end{align}

and $|\partial E|(u)$ is the metric slope defined by

\begin{equation}\label{metricDerDef}
|\partial E|(u) := \lim_{d(u,v) \to 0} \frac{(E(u) - E(v))_+}{d(u,v)} \in [0,\infty].
\end{equation}

\end{lemma}

\section{Statement of results}\label{sec:results}

Our main result is the convergence of the modified thresholding scheme to a weak notion of multiphase mean curvature flow. More precisely, given an initial partition $\{\Omega_1^0, . . . , \Omega_N^0\}$ of $[0,1)^d$ encoded by $\chi^0: [0,1)^d \to \{0,1\}^{N}$ such that $\sum_i \chi_i^0 = 1$, define $\chi^h: [0,1)^d \times \mathbf{R} \to \{0,1\}^N$ by setting

\begin{align}\label{defchih}
\begin{aligned}
&\chi^h(t,x) = \chi^0(x)\ \text{for}\ t < h,
\\ &\chi^h(t,x) = \chi^n(x)\ \text{for}\ t \in [nh, (n+1)h)\ \text{for}\ n \in \mathbf{N}.
\end{aligned}
\end{align}

If $\chi^0$ is a function of bounded variation, we denote by $\Sigma_{ij}^0 := \partial^*\Omega_i^0 \cap \partial^*\Omega_j^0$. Our main result is contained in the following theorem.
\begin{theorem}\label{mainThm}
Given $\chi^0 \in \mathcal{A}$ and such that $\nabla\chi^0$ is a bounded measure and a sequence $h\downarrow 0$; let $\chi^h$ be defined by (\ref{defchih}). Assume that there exists $\chi: [0,1)^d \times (0,T) \to [0,1]^N$ such that 

\begin{equation}\label{convergenceAssumpt}
\chi^h\rightharpoonup \chi\ \text{in}\ L^1([0,1)^d \times (0,T)).
\end{equation}
Then $\chi \in \{0,1\}^N$ almost everywhere, $\sum_{i}\chi_i = 1$ and $\chi \in L^1((0,T), BV([0,1)^d))^N$. 

If we assume that

\begin{equation}\label{energyAssumpt}
\limsup_{h \downarrow 0} \int_0^T E_h(\chi^h(t))dt \le \sum_{ij} \sigma_{ij}\int_0^T \Hm(\Sigma_{ij}(t))dt,
\end{equation}
then $\chi$ is a De Giorgi solution in the sense of Definition \ref{DeGiorgiSolution} below.
\end{theorem}
\begin{definition}\label{DeGiorgiSolution}
Given $\chi^0 \in \mathcal{A}$ and such that $\nabla\chi^0$ is a bounded measure, a map $\chi: [0,1)^d \times (0,T) \to \{0,1\}^N$ such that $\sum_{i}\chi_i = 1$ and $\chi \in L^1((0,T), BV([0,1)^d))^N$ is called a De Giorgi solution to the multiphase mean curvature flow with surface tensions $\sigma_{ij}$ and mobilities $\mu_{ij}$ provided the following three facts hold:
\begin{enumerate}
\item There exist $H_{ij} \in L^2(\Hm_{|\Sigma_{ij}}(dx)dt)$ which are mean curvatures in the weak sense, i.e., such that for any test vector field $\xi \in C^{\infty}_c([0,1)^d\times (0,T))^d$

\begin{align}\label{meanCUrvatureWeak}
&\sum_{i,j}\sigma_{ij}\int_{[0,1)^d\times (0,T)} (\nabla\cdot\xi - \nu_{ij}\cdot \nabla\xi \nu_{ij})\Hm_{|\Sigma_{ij}(t)}(dx)dt 
\\ &= -\sum_{i,j}\sigma_{ij}\int_{[0,1)^d\times (0,T)}H_{ij}\nu_{ij}\cdot\xi\Hm_{|\Sigma_{ij}(t)}(dx)dt.\nonumber
\end{align}
\item There exist normal velocities $V_{ij} \in L^2(\Hm_{|\Sigma_{ij}(t)}(dx)dt)$ with

\begin{align}
\begin{aligned}
\int_{[0,1)^d} \eta(t=0)\chi_i^0 dx &+ \int_{[0,1)^d\times (0,T)} \partial_t\eta\  \chi_i\ dxdt
\\ &+\sum_{k\neq i}\int_{[0,1)^d\times (0,T)} \eta V_{ik}\ \Hm_{|\Sigma_{ik}(t)}(dx)dt = 0\nonumber
\end{aligned}
\end{align}
for all $\eta \in C^{\infty}_c([0,1)^d\times [0,T))$.
\item De Giorgi's inequality is satisfied, i.e.\
\begin{align}\label{DEGIORGIFINAL}
\begin{aligned}
&\limsup_{\tau \downarrow 0}\frac{1}{\tau}\sum_{i,j}\sigma_{ij}\int_{(T-\tau, T)}\Hm(\Sigma_{ij}(t))dt
\\ &+ \frac{1}{2}\sum_{i,j}\int_{[0,1)^d\times (0,T)} \left(\frac{V_{ij}^2}{\mu_{ij}} + \sigma_{ij}^2\mu_{ij}H_{ij}^2 \right)\Hm_{|\Sigma_{ij}(t)}(dx)dt \le \sum_{i,j}\sigma_{ij}\Hm(\Sigma^0_{ij}).
\end{aligned}
\end{align}
\end{enumerate}
\end{definition}
\begin{remark}\label{relationSmoothSol}
Observe that inequality (\ref{DEGIORGIFINAL}) together with the definition of the weak mean curvatures gives a notion of weak solution for the multiphase mean curvature flow incorporating both the dynamics $V_{ij} = -\sigma_{ij}\mu_{ij}H_{ij}$ and the Herring angle condition at triple junctions. Indeed if $\chi: [0,1)^d \times (0,T) \to \{0,1\}^N$ with $\sum_{i} \chi_i(t) = 1$ is such that the sets $\Omega_i(t) = \{\chi_i(\cdot, t) = 1\}$ meet along smooth interfaces $\Sigma_{ij} := \partial \Omega_i \cap \partial \Omega_j$ which evolve smoothly and satisfy (\ref{meanCUrvatureWeak}), (\ref{DEGIORGIFINAL}) then
\begin{enumerate}
\item The Herring angle condition at triple junctions is satisfied. Indeed by the divergence theorem on surfaces (see Theorem 11.8 and Remark 11.42 in \cite{Maggi2012}) for any $\xi \in C^{\infty}_c(\times [0,1)^d)^d$ 
\begin{align*}
\int_{\Sigma_{ij}(t)} (\nabla\cdot\xi - \nu_{ij}\cdot\nabla\xi\nu_{ij})\ \Hm_{|\Sigma_{ij}(t)}(dx) = &-\int_{\Sigma_{ij}(t)} H_{ij}\nu_{ij}\Hm_{|\Sigma_{ij}(t)}(dx)
\\ & + \int_{\partial \Sigma_{ij}(t)} \xi \cdot \nu_{ij} \mathcal{H}^{d-2}(dx).
\end{align*}
Thus (\ref{meanCUrvatureWeak}) and $H_{ij} \in L^2(\Hm_{|\Sigma_{ij}(t)}(dx)dt)$ imply that 
\begin{align*}
&\sigma_{i_1i_2}\int_{\partial \Sigma_{i_1i_2}(t)} \xi \cdot \nu_{i_1i_2} \mathcal{H}^{d-2}(dx)
\\ & + \sigma_{i_2i_3}\int_{\partial \Sigma_{i_2i_3}(t)} \xi \cdot \nu_{i_2i_3} \mathcal{H}^{d-2}(dx)
\\ & + \sigma_{i_3i_1}\int_{\partial \Sigma_{i_3i_1}(t)} \xi \cdot \nu_{i_3i_1} \mathcal{H}^{d-2}(dx) = 0,
\end{align*}
which forces $\sigma_{i_1i_2}\nu_{i_1i_2} + \sigma_{i_2i_3}\nu_{i_2i_3} + \sigma_{i_3i_1}\nu_{i_3i_1} = 0$ at triple junctions.
\item We have $V_{ij} = -\sigma_{ij}\mu_{ij}H_{ij}$ on $\Sigma_{ij}(t)$. Indeed in the smooth case inequality (\ref{DEGIORGIFINAL}) reduces to
\begin{align*}
&\sum_{i,j}\sigma_{ij}\int_{(0, T)}\frac{d}{dt}\Hm(\Sigma_{ij}(t))dt
\\ &+ \frac{1}{2}\sum_{i,j}\int_{[0,1)^d\times (0,T)} \left(\frac{V_{ij}^2}{\mu_{ij}} + \sigma_{ij}^2\mu_{ij}H_{ij}^2 \right)\Hm_{|\Sigma_{ij}(t)}(dx)dt \le 0.
\end{align*}
If one recalls that $\frac{d}{dt}\Hm(\Sigma_{ij}(t)) = \int_{[0,1)^d} V_{ij}H_{ij} \Hm_{|\Sigma_{ij}(t)}(dx)$, after completing the square we arrive at 
\begin{align*}
\sum_{i,j}\int_{[0,1)^d\times (0,T)} \left(\frac{V_{ij}}{\sqrt{\mu_{ij}}} + \sigma_{ij}\mu_{ij}H_{ij} \right)^2\Hm_{|\Sigma_{ij}(t)}(dx)dt \le 0,
\end{align*}
which implies $V_{ij} = -\sigma_{ij}\mu_{ij}H_{ij}$.
\end{enumerate} 
\end{remark}
The following lemma establishes, next to a compactness statement, that our convergence can be localized in the space and time variables $x$ and $t$, but also in the variable $z$ appearing in the convolution.
\begin{lemma}\label{compactnessLemma}
\begin{enumerate}
\item  Let $\{\chi^h\}_{h\downarrow 0}$ be a sequence of $\{0,1\}^N$-valued functions on $(0,T) \times  [0,1)^d$ that satisfies

\begin{equation}\label{assumedBoundCom}
\limsup_{h\downarrow 0} \left( \esssup_{t \in (0,T)} E_h(\chi_h(t)) + \int_0^T \frac{1}{2h^2}d_h^2(\chi_h(t), \chi_h(t-h))dt\right) < \infty
\end{equation}
and that is piecewise constant in time in the sense of (\ref{defchih}). Such a sequence is compact in $L^1([0,1)^d \times (0,T))^N$ and any weak limit $\chi$ is such that $\chi \in L^1((0,T), \bv([0,1)^d))^N$ with 

\begin{equation}\label{liminfComp}
\sum_{i,j} \sigma_{ij}\int_0^T\Hm(\Sigma_{ij}(t))dt \le \liminf_{h \downarrow 0} \int_0^T E_h(\chi_h(t))dt.
\end{equation}

\item Assume that $u^h$ is a sequence of $[0,1]^N$-valued functions with $\sum_i u_i^h = 1$ such that (\ref{energyAssumpt}) holds (with $\chi^h$ replaced by $u^h$) and such that $u^h \to \chi$ in $L^1([0,1)^d\times (0,T))^N$ holds. Assume also that

\begin{equation}\label{energyBoundComp}
\limsup_{h\downarrow 0}\esssup_{t \in (0,T)} E_h(u^h(t)) < \infty.
\end{equation}
Then as measures on $\mathbf{R}^d \times [0,1)^d \times (0,T)$ we have the following weak convergences for any $i \neq j$

\begin{align}
&\begin{aligned}\label{provedConv}
\frac{K_{ij}(z)}{\sqrt{h}}u_i^h(x,t)&u_j^h(x-\sqrt{h}z,t)dxdtdz
\\ &\rightharpoonup K_{ij}(z)(\nu_{ij}(x,t)\cdot z)_+\Hm_{|\Sigma_{ij}(t)}(dx)dtdz.
\end{aligned}
\\ &\begin{aligned}\label{notprovedConv}
\frac{K_{ij}(z)}{\sqrt{h}}u_i^h(x-\sqrt{h}z,t)&u_j^h(x,t)dxdtdz
\\ &\rightharpoonup K_{ij}(z)(\nu_{ij}(x,t)\cdot z)_-\Hm_{|\Sigma_{ij}(t)}(dx)dtdz.
\end{aligned}
\end{align}
Here the convergence may be tested also with continuous functions which have polynomial growth in $z \in \mathbf{R}^d$.

\end{enumerate}
\end{lemma}
The next proposition is the main ingredient in the proof of Theorem \ref{mainThm}. It establishes the sharp lower bound on the distance-term.
\begin{proposition}\label{propVelocity}

Suppose that (\ref{convergenceAssumpt}) and the conclusion of Lemma \ref{compactnessLemma} (2) hold. Assume also that the left hand side of (\ref{sharpVelocityBound}) is finite. Then for every $1 \le k \le N$ there exists $V_k \in L^2(|\nabla\chi_k|dt)$ such that

\begin{equation}\label{normalVelocity}
\partial_t\chi_k = V_k|\nabla\chi_k|dt
\end{equation}
in the sense of distributions. Given $i \neq j$, it holds that $V_{i}(x,t) = - V_{j}(x,t)$ on $\Sigma_{ij}(t)$ and if we define $V_{ij}(x,t) := V_i(x, t) \nu_{ij}(x,t)|_{\Sigma_{ij}(t)}$ then we have

\begin{equation}\label{sharpVelocityBound}
\liminf_{h \downarrow 0} \int_{0}^T \frac{1}{h^2}d_h^2(\chi^h(t), \chi^h(t-h))dt \ge \sum_{i,j}\frac{1}{\mu_{ij}} \int_0^T\int_{\Sigma_{ij}(t)} |V_{ij}(x,t)|^2\Hm(dx)dt.
\end{equation}

\end{proposition}
The final ingredient is the analogous sharp lower bound for the metric slope.
\begin{proposition}\label{meanCurv}
Suppose that the conclusion of Lemma \ref{compactnessLemma} (2) holds and that (\ref{convergenceAssumpt}) holds with $\chi^h$ replaced by $u^h$. Then for any $i \neq j$ there exists a mean curvature $H_{ij} \in L^2(\Hm_{\Sigma_{ij}(t)}(dx)dt)$ in the sense of (\ref{meanCUrvatureWeak}). Moreover the following inequality is true:

\begin{equation}\label{meanCurvBound}
\liminf_{h\downarrow 0} \int_0^T |\partial E_h|^2(u_h(t)) dt \ge \sum_{i,j} \mu_{ij}\sigma_{ij}^2\int_0^T\int_{\Sigma_{ij}(t)}|H_{ij}(x,t)|^2\Hm(dx)dt.\nonumber
\end{equation}

\end{proposition}

\section{Construction of suitable partitions of unity}\label{consPOU}

In the sequel we will frequently want to localize on one of the interfaces. To do so, we need to construct a suitable family of balls on wich the behavior of the flow is split into two majority phases and several minority phases. Hereafter we will ignore the time variable and consider a map $\chi : [0, 1)^d \to \{0,1\}^N$ such that $\chi \in BV([0,1)^d, \mathbf{R}^N)$, $\sum_k \chi_k = 1$. Given $1 \le i < j \le N$ we denote by $\partial^*\Omega_i$ the reduced boundary of the set $\{\chi_i = 1\}$ and by $\Sigma_{ij} = \partial^*\Omega_i \cap \partial^*\Omega_j$ the interface between phase $i$ and phase $j$. Given a real number $r > 0$ and a natural number $n \in \mathbf{N}$ we define

\begin{equation}\label{ballsGrid}
\mathcal{F}_n^r:=\left\{ B(x, nr\sqrt{d}):\ x \in r\mathbf{Z}^d\right\}
\end{equation}
where the balls appearing in the definition are intended to be open. Observe that for any $n \ge 2$ and any $r > 0$ the collection of balls in $\mathcal{F}_n^r$ is a covering of $\mathbf{R}^d$ with the property that any point $x \in \mathbf{R}^d$ lies in at most $c(n,d)$ distinct balls belonging to $\mathcal{F}^r_n$, where $0 < c(n,d) \le (2n)^d$ is a constant that depends on $n, d$ but not on $r$. Given numbers $1 \le l \neq p \le N$ we define

\begin{equation}\label{goodBalls}
\mathcal{E}^r := \left\{ B \in \mathcal{F}^r_2:\ B \cap \Sigma_{lp} \neq \emptyset,\ \frac{\Hm(\Sigma_{ij} \cap 2B)}{\omega_{d-1}(4r)^{d-1}} \le \frac{1}{2^d},\ \{i,j\} \neq \{l,p\} \right\}.
\end{equation}
Here $2B$ denotes the ball with center given by the center of $B$ and twice its radius.
Given $l,p$ as above, denote by $\{B_{m}^r\}_{m\in\mathbf{N}}$ an enumeration of $\mathcal{E}^r$ and by $\{\rho_{m}\}_{m\in\mathbf{N}}$ a smooth partition of unity subordinate to $\{B_{m}^r\}_{m\in\mathbf{N}}$.

\begin{lemma}\label{propertiesPOU}

Fix $1 \le l \neq p \le N$. With the above construction the following two properties hold.

\begin{enumerate}

\item For any $1 \le i \neq j \le N$, $\{i,j\} \neq \{l,p\}$ and any $\eta \in L^1(\Hm_{|\Sigma_{ij}})$

\begin{equation}\label{errorBound}
\lim_{r\downarrow 0} \sum_{m \in \mathbf{N}} \int_{B_{m}^r} \eta \Hm_{|\Sigma_{ij}}(dx) = 0.
\end{equation}

\item For any $\eta \in L^1(\Hm_{|\Sigma_{lp}})$

\begin{equation}\label{approximationInterface}
\lim_{r\downarrow 0} \sum_{m \in \mathbf{N}} \int \rho_{m} \eta \Hm_{|\Sigma_{lp}}(dx) = \int \eta \Hm_{|\Sigma_{lp}}(dx).
\end{equation}

\end{enumerate}

\end{lemma}

\section{Proofs}\label{sec:proofs}

\begin{proof}[Proof of Theorem \ref{mainThm}]
By Lemma \ref{minmovInt}, we can apply Lemma \ref{afpLemma} on the metric space $(\mathcal{M}, d_h)$ so that we get inequality (\ref{DGineq}) with $(E,d,\chi,u) = (E_h, d_h, \chi^h, u^h)$. Our first observation is that

\begin{equation}\label{RHS}
\lim_{h\downarrow 0} E_h(\chi^0) = \sum_{i,j}\sigma_{ij}\Hm(\Sigma_{ij}^0),
\end{equation}
which follows from the consistency, cf.\ Lemma \ref{consistency} in the Appendix. Inequality (\ref{DGineq}) then yields that the sequence $\chi^h$ satisfies (\ref{assumedBoundCom}), so that Lemma \ref{compactnessLemma} (1) applies to get that $\chi \in L^1((0,T), BV([0,1)^d))^N$, $\chi \in \{0,1\}^N$ a.e., $\sum_{i} \chi_i = 1$ and, after extracting a subsequence, $\chi^h \to \chi$ in $L^2([0,1)^d \times (0,T))^N$. We claim that this implies $u^h \to \chi$ in $L^2([0,1)^d \times (0,T))^N$. To see this, observe that (\ref{varDist}) implies

\begin{align}\label{closeUX}
\begin{aligned}
hE_h(\chi^0) &\ge -\int_0^T E_h(u_h(t) - \chi_h(t))dt
\\ &\ge C \frac{1}{\sqrt{h}}\sum_{i=1}^N\left(\int |G_{\gamma}^{h/2}*(u_i^h - \chi_i^h)|^2 dxdt + \int|G_{\beta}^{h/2}*(u_i^h - \chi_i^h)|^2 dxdt \right)
\end{aligned}
\end{align}
where $C$ is a constant which depends on $N, \mathbb{A}, \mathbb{B}$ but not on $h$ and comes from the fact that all norms on $(1, . . . , 1)^{\perp}$ are comparable. Inequality (\ref{closeUX}) clearly implies that $K^h*u_h - K^h*\chi_h$ converges to zero in $L^2$. Observe that inequality (\ref{varEnergy}) in particular yields (\ref{energyBoundComp}) with $\chi^h$ replaced by $u^h$. Recalling (\ref{closeL2K}) in the Appendix, we learn that $u^h - \chi^h$ converges to zero in $L^2$.  This implies that we can apply Lemma \ref{compactnessLemma} (2) both to the sequence $u^h$ and the sequence $\chi^h$. In particular, we may apply Proposition \ref{propVelocity} for $\chi^h$ and Proposition \ref{meanCurv} for $u^h$. Now the proof follows the same strategy as the one in the two-phase case in \cite{Laux2019}. For the sake of completeness, we sketch the argument here. First of all, Lemma \ref{afpLemma} gives inequality (\ref{DGineq}) for $(E_h, d_h, \chi^h, u^h)$, namely for $n \in \mathbf{N}$
\begin{align}\label{trickrho}
\rho(nh) \le E_h(\chi^0),
\end{align}
where we set $\rho(t) = E_h(\chi^h(t)) + \frac{1}{2}\int_0^t \left( \frac{1}{h^2} d_h^2(\chi^h(s+h), \chi^h(s)) + |\partial E_h(u^h(s))|^s\right) ds$. Multiplying (\ref{trickrho}) by $\eta(nh) - \eta((n+1)h)$ for some non-increasing function $\eta \in C_c([0,T))$ we get $-\int \frac{d\eta}{dt} \rho dt \le (\eta(0) + h \sup \left| \frac{d\eta}{dt} \right|)E_h(\chi^0)$. As test function $\eta$, we now choose $\eta(t) = \max\{\min\{\frac{T-t}{\tau}, 1\}, 0\}$ and obtain
\begin{align}\label{beforePassing}
& \frac{1}{\tau}\int_{T-\tau}^T E_h(\chi^h(t))dt 
\\ &+ \frac{1}{2}\int_0^{T-\tau} \left(\frac{1}{h^2}d_h^2(\chi^h(t), \chi^h(t-h)) + |\partial E_h(u^h(t))|^2\right)dt \le (1 + \frac{h}{\tau})E_h(\chi^0).\nonumber
\end{align}
Now it remains to pass to the limit as $h \downarrow 0$: to get (\ref{DEGIORGIFINAL}) from inequality (\ref{beforePassing}) one uses the lower semicontinuity (\ref{liminfComp}) for the first left hand side term, the sharp bound (\ref{sharpVelocityBound}) for the second left hand side term, the bound (\ref{meanCurvBound}) for the last left hand side term and finally one uses the consistency Lemma \ref{consistency} in the Appendix to treat the right hand side term. To get (\ref{DEGIORGIFINAL}) it remains to pass to the limit in $\tau \downarrow 0$.

\end{proof}

\begin{proof}[Proof of Lemma \ref{compactnessLemma}]
Argument for (1). For the compactness, the arguments in \cite{Laux2019} adapt to this setting with minor changes. The first observation is that, by inequality (\ref{closeL2K}) in the Appendix, one needs to prove compactness in $L^2([0,1)^d \times (0,T))^N$ of $\{K^h * \chi^h\}_{h\downarrow 0}$. For this, one just needs a modulus of continuity in time. I.e.\ it is sufficient to prove that there exists a constant $C > 0$ independent of $h$ such that $I_h(s) \le C\sqrt{s}$, where

\begin{equation}
I_h(s) = \int_{(s,T) \times [0,1)^d} |\chi_h(x,t) - \chi_h(x, t-s)|^2dxdt.\nonumber
\end{equation}

This is can be done applying word by word the argument in \cite{Laux2019} once we show the following: for any pair $\chi, \chi' \in \mathcal{A}$ of admissible functions, we have

\begin{equation}\label{ineqForMC}
\int |\chi-\chi'|dx \le \frac{C}{\sqrt{h}}d_h^2(\chi, \chi') + C\sqrt{h}\left(E_h(\chi) + E_h(\chi')\right).
\end{equation}
Here the constant $C$ depends on $N, \mathbb{A}, \mathbb{B}$ but not on $h$.

To prove (\ref{ineqForMC}) we proceed as follows: let $\mathbb{S} \in \mathbf{R}^{N\times N}$ be a symmetric matrix which is positive definite on $(1, . . . , 1)^{\perp}$. Since any two norms on a finite dimensional space are comparable, there exists a constant $C>0$ depending on $\mathbb{S}$ and $N$ such that

\begin{equation}
|\chi - \chi'| \le |\chi - \chi'|^2 \le C|\chi - \chi'|^2_{\mathbb{S}}
\end{equation}
where $|\cdot|_{\mathbb{S}}$ denotes the norm induced by $\mathbb{S}$. For a function $u \in \mathcal{M}$ write $(\tilde{K}^h *) u_h$ for the function

\begin{equation}
\left((\tilde{K}^h *) u_h\right)_i = \sum_{j \neq i} K_{ij}^h * u_h^j.\nonumber
\end{equation}
Then we calculate

\begin{align}\label{smugglingS}
\begin{aligned}
|\chi-\chi'|^2_{\mathbb{S}} = -(\chi-\chi')\cdot (\tilde{K}^h*)(\chi-\chi') + (\chi-\chi')(\mathbb{S} + (\tilde{K}^h*))(\chi-\chi').
\end{aligned}
\end{align}
Select $\mathbb{S} = (s_{ij})$ where $s_{ij} = -\int K_{ij}(z) dz$. Then, by our assumption (\ref{posDefAss}) $\mathbb{S}$ is positive definite on $(1, . . . , 1)^{\perp}$ and after integration on $[0,1)^d$ identity (\ref{smugglingS}) becomes

\begin{equation}
\int |\chi-\chi'|^2_{\mathbb{S}}dx = \frac{1}{2\sqrt{h}}d_h^2(\chi, \chi') + \int (\chi-\chi')(\mathbb{S} + (\tilde{K}^h*))(\chi-\chi') dx.\nonumber
\end{equation}
We now proceed to estimate the integral on the right hand side. By the choice of $\mathbb{S}$ and Jensen's inequality we have

\begin{align}
&\int (\chi-\chi')(\mathbb{S} + (\tilde{K}^h*))(\chi-\chi') dx \le C\int |\mathbb{S} + (\tilde{K}^h*))(\chi-\chi')| dx
\\ & \le  C\sum_{i,j} \int K^h_{ij}(z)|(\chi_j-\chi_j')(x-z) - (\chi_j-\chi_j')(x)|dxdz.\nonumber
\end{align}
Using the triangle inequality and (\ref{ineqAB}) in the Appendix we can estimate the right hand side to obtain the following inequality

\begin{align}
&\begin{aligned}
\int (\chi-\chi')(\mathbb{S} + (\tilde{K}^h*))(\chi-\chi') dx
\end{aligned} 
\\ &\begin{aligned}
\le C\sum_{i,j} &\left(\sum_{k \neq j}\int K^h_{ij}(z)\chi_j(x-z)\chi_k(x)dxdz\right.
\\ &\left. +\sum_{k \neq j}\int K^h_{ij}(z)\chi_j(x)\chi_k(x-z)dxdz\right.
\\ &\left. +\sum_{k \neq j}\int K^h_{ij}(z)\chi'_j(x-z)\chi'_k(x)dxdz\right. 
\\ &\left. + \sum_{k \neq j}\int K^h_{ij}(z)\chi'_j(x)\chi'_k(x-z)dxdz\right).\nonumber
\end{aligned}
\end{align}
Observing that there is a constant $C>0$ such that $K_{ij} \le CK_{jk}$ we conclude that

\begin{equation}
\int (\chi-\chi')(\mathbb{S} + (\tilde{K}^h*))(\chi-\chi') dx \le C \sqrt{h}\left(E_h(\chi) + E_{h}(\chi') \right).\nonumber
\end{equation}
This proves (\ref{ineqForMC}) and closes the argument for the compactness.

We also have to prove (\ref{liminfComp}), but this follows from (\ref{provedConv}) with $u^h$ replaced by $\chi^h$ once we have shown that the limit $\chi$ is such that $|\nabla \chi|$ is a bounded measure, equiintegrable in time. This can be done with an argument similar to the one used in \cite{Laux2019} for the two-phase case. Observe that this only requires the weaker assumption (\ref{energyBoundComp}).

Argument for (2). As mentioned in the previous paragraph, we already know that the limit $\chi$ is such that $|\nabla \chi|$ is a bounded measure, equiintegrable in time. We will prove (\ref{provedConv}). Then (\ref{notprovedConv}) easily follows by recalling that $\nu_{ij} = - \nu_{ji}$.
A standard argument (to be found in \cite{Laux2019}) which relies on the exponential decay of the kernel yields the fact that we can test convergences (\ref{provedConv}) with functions with at most polynomial growth in $z$ provided we already have the result for bounded and continuous test functions, thus we focus on this case. 

Let $\xi \in C_b(\mathbf{R}^d \times [0,1)^d \times (0,T))$ be a bounded and continuous function. To show (\ref{provedConv}) we aim at showing that

\begin{align}
\begin{aligned}\label{convToSplit}
\lim_{h\downarrow 0}\int\xi(z,x,t)\frac{K_{ij}(z)}{\sqrt{h}}&u_i^h(x,t)u_j^h(x-\sqrt{h}z,t)dxdtdz
\\ &= \int\xi(z,x,t)K_{ij}(z)(\nu_{ij}(x,t)\cdot z)_+ \Hm_{\Sigma_{ij}(t)}(dx)dtdz.
\end{aligned}
\end{align}

Upon splitting $\xi$ into the positive and the negative part, by linearity we may assume that $0 \le \xi \le 1$. We can split (\ref{convToSplit}) into the local lower bound

\begin{align}
\begin{aligned}\label{convToSplitted}
\liminf_{h\downarrow 0}\int\xi(z,x,t)\frac{K_{ij}(z)}{\sqrt{h}}&u_i^h(x,t)u_j^h(x-\sqrt{h}z,t)dzdxdt
\\ &\ge \int\xi(z,x,t)K_{ij}(z)(\nu_{ij}(x,t)\cdot z)_+ \Hm_{\Sigma_{ij}(t)}(dx)dtdz.
\end{aligned}
\end{align}
and the global upper bound

\begin{align}
\begin{aligned}\label{convToSplittedUp}
\liminf_{h\downarrow 0}\int\frac{K_{ij}(z)}{\sqrt{h}}&u_i^h(x,t)u_j^h(x-\sqrt{h}z,t)dzdxdt
\\ &\le \int K_{ij}(z)(\nu_{ij}(x,t)\cdot z)_+ \Hm_{\Sigma_{ij}(t)}(dx)dtdz.
\end{aligned}
\end{align}
Indeed we can recover the limsup inequality in (\ref{convToSplit}) by splitting $\xi = 1 - (1 - \xi)$ and applying the local lower bound (\ref{convToSplitted}) to $1-\xi$. 

We first concentrate on the local lower bounds in the case where $u^h = \chi$, namely we will show

\begin{align}
\begin{aligned}\label{no_h}
\liminf_{h\downarrow 0}\int\xi(z,x,t)\frac{K_{ij}(z)}{\sqrt{h}}&\chi_i(x,t)\chi_j(x-\sqrt{h}z,t)dzdxdt
\\ &\ge \int\xi(z,x,t)K_{ij}(z)(\nu_{ij}(x,t)\cdot z)_+ \Hm_{\Sigma_{ij}(t)}(dx)dtdz.
\end{aligned}
\end{align}

By Fatou's lemma the claim is reduced to showing that for a.e.\ point $t$ in time and every $z \in \mathbf{R}^d$

\begin{align}
\liminf_{h\downarrow 0}\int\xi(z,x,t)\frac{K_{ij}(z)}{\sqrt{h}}&\chi^i(x,t)\chi^j(x-\sqrt{h}z,t)dx
\\ &\ge \int\xi(z,x,t)K_{ij}(z)(\nu_{ij}(x,t)\cdot z)_+ \Hm_{\Sigma_{ij}(t)}(dx).\nonumber
\end{align}

Fix a point $t$ such that $\chi(\cdot, t) \in \bv([0,1)^d, \{0,1\}^N)$ and any $z \in \mathbf{R}^d$. In the sequel, we will drop those variables, so $\chi(x) = \chi(x,t)$, $\xi(x) = \xi(z,x,t)$. By approximation we may assume that $\xi \in C^{\infty}([0,1)^d)$. Let $\rho_{m}$ be a partition of unity obtained by applying the construction of Section \ref{consPOU} to the function $\chi(x)$ on the interface $\Sigma_{ij}$. Then by Lemma \ref{propertiesPOU} we have

\begin{align}
&\begin{aligned}
&\int \xi(x) (\nu_{ij}(x) \cdot z)_+\Hm_{|\Sigma_{ij}}(dx)
\\ & = \lim_{r\downarrow 0} \left( \sum_{m \in \mathbf{N}} \int \rho_{mij}(x)\xi(x) (\nu_{ij}(x) \cdot z)_+\Hm_{|\Sigma_{ij}}(dx)\right)
\end{aligned}\nonumber
\\ & \begin{aligned}\nonumber
= \lim_{r\downarrow 0} \sum_{m \in \mathbf{N}}  &\left( \int \rho_{mij}(x)\xi(x) (\nu_{i}(x) \cdot z)_+\Hm_{|\partial^*\Omega_i}(dx)\right.
\\ &\left. - \sum_{k \neq i,j}\int \rho_{mij}(x)\xi(x) (\nu_{ij}(x) \cdot z)_+\Hm_{|\Sigma_{ik}}(dx) \right)
\end{aligned}\nonumber
\\ & \begin{aligned}\label{limToEst}
= \lim_{r\downarrow 0} \sum_{m \in \mathbf{N}} \int \rho_{mij}(x)\xi(x) (\nu_{i}(x) \cdot z)_+\Hm_{|\partial^*\Omega_i}(dx).
\end{aligned}
\end{align}
We now focus on estimating the argument of the last limit. Observe that $(\nu_{ij}(x) \cdot z)_+\Hm_{|\partial^*\Omega_i}(dx) = (\partial_z \chi_i)_+$, thus by definition of positive part of a measure, given $\epsilon > 0$ we can select, for any $m \in \mathbf{N}$, a function $\tilde{\xi}_m \in C^1_c(B_{m})$ such that $0 \le \tilde{\xi}_m \le 1$ and such that

\begin{equation}\label{choiceXi}
\int \rho_{mij}\xi\tilde{\xi}_m \partial_z\chi_i + 2^{-m}\epsilon \ge \int\rho_{mij}\xi(\nu_i \cdot z)_+\Hm_{|\partial^*\Omega_i}(dx).
\end{equation}
Let $\eta_m := \rho_{mij}\xi\tilde{\xi}_m \in C^1_c(B_m)$, then

\begin{align}
&\begin{aligned}
\int\eta_m\partial_z\chi_i = -\int\partial_z\eta_m\chi_i dx
\end{aligned}\nonumber
\\ &\begin{aligned}
= \lim_{h\downarrow 0} \int \frac{\eta_m(x+\sqrt{h}z) - \eta(x)}{\sqrt{h}} \chi_i(x)dx
\end{aligned}\nonumber
\\ &\begin{aligned}
= \lim_{h\downarrow 0} \int \eta_m(x) \frac{\chi_i(x) - \chi_i(x-\sqrt{h}z)}{\sqrt{h}}dx
\end{aligned}\nonumber
\\ &\begin{aligned}
\le \liminf_{h\downarrow 0} \sum_{k\neq i}\int \eta_m(x) \frac{\chi_i(x)\chi_k(x-\sqrt{h}z)}{\sqrt{h}}dx
\end{aligned}\nonumber
\\ &\begin{aligned}
\le \liminf_{h\downarrow 0} &\int \eta_m(x) \frac{\chi_i(x)\chi_j(x-\sqrt{h}z)}{\sqrt{h}}dx 
\\ &+ \limsup_{h\downarrow 0}  \sum_{k\neq i,j}\int \eta_m(x) \frac{\chi_i(x)\chi_k(x-\sqrt{h}z)}{\sqrt{h}}dx
\end{aligned}\nonumber
\\ &\begin{aligned}
\le \liminf_{h\downarrow 0} &\int \eta_m(x) \frac{\chi_i(x)\chi_j(x-\sqrt{h}z)}{\sqrt{h}}dx 
\\ &+ \sum_{k\neq i,j}\limsup_{h\downarrow 0}\int \eta_m(x) \frac{\chi_i(x)\chi_k(x-\sqrt{h}z)}{\sqrt{h}}dx.
\end{aligned}\nonumber
\end{align}
Observe that for each $m \in \mathbf{N}$, using also the consistency Lemma \ref{consistency}

\begin{align*}
&\begin{aligned}
\limsup_{h\downarrow 0}\int \eta_m(x) \frac{\chi_i(x)\chi_k(x-\sqrt{h}z)}{\sqrt{h}}dx
\end{aligned}
\\ &\begin{aligned}
\le \limsup_{h\downarrow 0}\int \eta_m(x) \frac{\chi_i(x)\chi_k(x-\sqrt{h}z)+\chi_i(x-\sqrt{h}z)\chi_j(x)}{\sqrt{h}}dx
\end{aligned}
\\ &\begin{aligned}
= \int \eta_m(x) |\nu_{ik}(x) \cdot z| \Hm_{|\Sigma_{ik}}(dx)
\end{aligned}
\\ &\begin{aligned}
\le |z|\Hm(B_{mij}^r \cap \Sigma_{ik})
\end{aligned}
\end{align*}
Thus we obtain

\begin{align}
\begin{aligned}
\int\eta_m\partial_z\chi_i \le \liminf_{h\downarrow 0} &\int \eta_m(x) \frac{\chi_i(x)\chi_j(x-\sqrt{h}z)}{\sqrt{h}}dx 
\\ &+ \sum_{k\neq i,j}|z|\Hm(B_{mij}^r \cap \Sigma_{ik})\nonumber
\end{aligned}
\end{align}

Inserting back into (\ref{limToEst}), recalling also Lemma \ref{propertiesPOU} and the inequality (\ref{choiceXi}), using Fatou's lemma, the fact that $\rho_{mij}$ is a partition of unity and that $0\le \tilde{\xi}_m \le 1$ we obtain that

\begin{align}
\int \xi(x) (\nu_{ij}(x) \cdot z)_+\Hm_{|\Sigma_{ij}}(dx) \le \liminf_{h\downarrow 0} &\int \xi(x) \frac{\chi_i(x)\chi_j(x-\sqrt{h}z)}{\sqrt{h}}dx + \epsilon\nonumber
\end{align}
and (\ref{no_h}) follows letting $\epsilon$ go to zero. To derive inequality (\ref{convToSplitted}) we just apply Lemma \ref{noHtoH} in the Appendix.

To get the upper bound (\ref{convToSplittedUp}) we argue as follows. First of all recall Assumption  (\ref{energyAssumpt}) which says

\begin{equation}\label{energyCOnv}
\int_0^TE_h(u^h(t))dt \to \int_0^T E(\chi(t))dt.
\end{equation}
Now, if we define 

\begin{equation}
e^{ij}_h(u^h) = \frac{1}{\sqrt{h}}\int_0^T\int u_i^h(t) K_{ij}^h * u_j^h(t) dxdt\nonumber
\end{equation}
we have that by (\ref{convToSplitted}) $\liminf_{h\downarrow 0} e_h^{ij}(u_h) \ge e^{ij}(\chi)$, where $e^{ij}(\chi)$ is defined in the obvious way. Assume that there exists a pair $i,j$ such that $\limsup_{h\downarrow 0} e_h^{ij}(u^h) > e^{ij}(\chi)$, then

\begin{align}
\begin{aligned}
\int_0^T E(\chi(t))dt &= \lim_{h\downarrow 0}\int_0^TE_h(u^h(t))dt
\\ & = \limsup_{h\downarrow 0} \int_0^TE_h(u^h(t))dt
\\ & \ge \sum_{(l,p) \neq (i,j)} \liminf_{h\downarrow 0} e_h^{lp}(u^h) + \limsup_{h\downarrow 0} e_h^{ij}(u^h)
\\ & > \int_0^T E(\chi(t))dt\nonumber
\end{aligned}
\end{align}
which is a contradiction. Thus we have proved (\ref{convToSplittedUp}).
\end{proof}

\begin{proof}[Proof of Proposition \ref{propVelocity}]
Since we assume that the left hand side of (\ref{sharpVelocityBound}) is finite, in view of (\ref{distanceDef}), upon passing to a subsequence we may assume that, in the sense of distributions, the limit

\begin{equation}\label{dissipationMeasure}
\lim_{h\downarrow 0} \frac{1}{h\sqrt{h}} \left( \left| G^{h/2}_{\gamma} * (\chi - \chi( \cdot - h))\right|_{\mathbb{A}}^2 + \left| G^{h/2}_{\beta}*(\chi - \chi( \cdot - h))\right|_{\mathbb{B}}^2 \right) = \omega
\end{equation}
exists as a finite positive measure on $[0,1)^d \times (0,T)$. Here we indicated with $\chi_l^h( \cdot - h)$ the time shift of function $\chi_l^h$. We denote by $\tau$ a small fraction of the characteristic spatial scale, namely $\tau = \alpha \sqrt{h}$ for some $\alpha > 0$, which we think as a small number. Given $1 \le l \le N$ we define

\begin{equation}\label{deltaDefinition}
\delta\chi^h_l := \chi^h_l - \chi^h_l( \cdot - \tau).
\end{equation}

We divide the proof into two parts: first we show that the normal velocities exist, and afterwards we prove the sharp bound. But first, let us state two distributional inequalities that will be used later. Namely

\begin{itemize}

\item In a distributional sense it holds that

\begin{equation}\label{DissMeasBound}
\limsup_{h\downarrow 0} -\frac{1}{\sqrt{h}}\sum_{i\neq j}\delta\chi_i K_{ij}^h *\delta\chi_j \le \alpha^2\omega.
\end{equation}

\item There exists a constant $C> 0$ such that for any $1 \le i \le N$ and any $\theta \in \{\gamma, \beta\}$ in a distributional sense it holds that

\begin{equation}\label{DissMeasBoundOneInt}
\limsup_{h\downarrow 0} \frac{1}{\sqrt{h}}(\chi_i - \chi_i(\cdot-\tau))G_{\theta}^h * (\chi_i - \chi_i(\cdot-\tau)) \le C\alpha^2\omega.
\end{equation}

\end{itemize}

We observe that it suffices to prove (\ref{DissMeasBound}), then (\ref{DissMeasBoundOneInt}) follows immediately. Indeed recall that $\mathbb{A}$ and $\mathbb{B}$ are positive definite on $(1, . . . , 1)^{\perp}$. In particular there exists a constant $C>0$ such that for any $v \in (1, . . . , 1)^{\perp}$ one has $|v|^2_{\mathbb{A}} + |v|_{\mathbb{B}}^2 \ge C|v|^2 \ge Cv_i^2$ for any $i \in \{1, . . . , N\}$. Applying this to the vector $v = G_\theta^{h/2}*\delta\chi_i$ one gets 

\begin{equation}
|G_\theta^{h/2}*\delta\chi_i|^2 \le \frac{1}{C} |G_{\theta}^{h/2} *\delta\chi|_{\mathbb{A}}^2 + |G_{\theta}^{h/2} * \delta\chi|_{\mathbb{B}}^2.
\end{equation}
The claim then follows from the definition of $\omega$, (\ref{DissMeasBound}), the symmetry (\ref{symmetryHK}) and the semigroup property (\ref{semigroupHK}). Indeed it is sufficient to check that, in the sense of distributions

\begin{equation}\label{HadTrick}
\lim_{h\downarrow 0} \frac{1}{\sqrt{h}}\sum_{i \neq j} \delta\chi_i K_{ij}^h*\delta\chi_j +\frac{1}{\sqrt{h}}\left(|G_{\gamma}^{h/2} *\delta\chi|_{\mathbb{A}}^2 + |G_{\beta}^{h/2} * \delta\chi|_{\mathbb{B}}^2 \right) = 0.
\end{equation}

To this aim, pick a test function $\eta \in C^{\infty}_c([0,1)^d\times (0,T))$. Spelling out the definition of the norms $|\cdot|_{\mathbb{A}}$ and $|\cdot|_{\mathbb{B}}$, the claim is proved once we show that

\begin{align}\label{eqProvHadT}
\lim_{h\downarrow 0} \frac{1}{\sqrt{h}}\sum_{i,j} a_{ij} \int \xi(\delta\chi_i G_{\gamma}^h*\delta\chi_j - G_{\gamma}^{h/2}*\delta\chi_i G_{\gamma}^{h/2}*\delta\chi_j)dxdt = 0,
\end{align}
and the same claim with $a_{ij}$, $\gamma$ replaced by $b_{ij}, \beta$ respectively. 

We concentrate on (\ref{eqProvHadT}). Clearly, we are done once we show that for any $i \neq j$ 

\begin{align}\label{toshowHT}
\lim_{h\downarrow 0} \frac{1}{\sqrt{h}}\int \xi(\delta\chi_i G_{\gamma}^h*\delta\chi_j - G_{\gamma}^{h/2}*\delta\chi_i G_{\gamma}^{h/2}*\delta\chi_j)dxdt = 0.
\end{align}
To show this, using the semigroup property (\ref{semigroupHK}) we rewrite the argument of the limit as

\begin{equation}
-\frac{1}{\sqrt{h}}\int [\xi, G_{\gamma}^{h/2}*](\delta\chi_i)G_{\gamma}^{h/2}*\delta\chi_j dxdt,
\end{equation}
and we observe that by the boundedness of the measures $\frac{1}{\sqrt{h}}|G_{\gamma}^{h/2}*\delta\chi|^2_{\mathbb{A}}$ it suffices to show

\begin{equation}\label{limitToZero}
\lim_{h\downarrow 0} \frac{1}{\sqrt{h}}\int|[\xi, G_{\gamma}^{h/2}*](\delta\chi_i)|^2dxdt = 0.
\end{equation}
To prove this, spelling out the integrand, using the Cauchy-Schwarz inequality and recalling the scaling (\ref{scalingHK}) we observe that

\begin{align}
&\begin{aligned}
\int|[\xi, G_{\gamma}^{h/2}*](\delta\chi_i)|^2dxdt 
\end{aligned}\nonumber
\\ &\begin{aligned}\label{boundOh}
&\le \int\left(\int|\xi(x,t) - \xi(x-z,t)|^2G_{\gamma}^{h/2}(z)dz\right)G_{\gamma}^{h/2}*|\delta\chi_i(x,t)|^2dxdt
\\ &\le \frac{h}{2} \sup |\nabla \xi|^2\int G_{\gamma}(z)|z|^2dz \int_0^T \int |\delta\chi_i(x,t)|^2dxdt.
\end{aligned}
\end{align}
Observe that by the compactness of $\chi^h$ in $L^2([0,1)^d\times (0,T))$, (\ref{boundOh}) is of order $h$, thus (\ref{limitToZero}) indeed holds true.

The proof of (\ref{DissMeasBound}) is essentially already contained in the paper \cite{Laux2019}. For the convenience of the reader we sketch the main ideas here. One reduces the claim to proving the following facts. 

\begin{equation}\label{subst1}
\lim_{h\downarrow 0} \frac{1}{\sqrt{h}}\sum_{ij}\delta\chi_i K^h_{ij}*\delta\chi_j - \frac{1}{\sqrt{h}}\left(\left| G_{\gamma } ^{h/2}* \delta\chi \right|_{\mathbb{A}}^2 + \left| G_{\beta}^{h/2}*\delta\chi\right|_{\mathbb{B}}^2\right) = 0.
\end{equation}

\begin{equation}\label{subst2}
\limsup_{h \downarrow 0} \frac{1}{\sqrt{h}}\left| G_{\gamma}^{h/2} * \delta\chi \right|_{\mathbb{A}}^2 - \alpha^2 \frac{1}{h\sqrt{h}}\left| G_{\gamma}^{h/2} * (\chi - \chi( \cdot - h))\right|_{\mathbb{A}}^2 \le 0.
\end{equation}

\begin{equation}\label{subst3}
\limsup_{h \downarrow 0} \frac{1}{\sqrt{h}}\left| G_{\beta}^{h/2}*\delta\chi\right|_{\mathbb{B}}^2 - \alpha^2 \frac{1}{h\sqrt{h}} \left| G_{\beta}^{h/2}*(\chi - \chi( \cdot - h))\right|_{\mathbb{B}}^2 \le 0.
\end{equation}

Claim (\ref{subst1}) was proved in the previous paragraph, while (\ref{subst2}) and (\ref{subst3}) are consequences of Jensen's inequality in the time variable for the convex functions $|\cdot|_{\mathbb{A}}^2$ and $|\cdot|_{\mathbb{B}}^2$ respectively. More precisely, assume without loss of generality that $\tau = Nh$ for some $N \in \mathbf{N}$, then by a telescoping argument and Jensen's inequality for $|\cdot|_{\mathbb{A}}^2$ we get

\begin{align}
\begin{aligned}
&\frac{1}{\sqrt{h}}|G_{\gamma}^{h/2}*\delta\chi|_{\mathbb{A}}^2
\\ &\le N\sum_{n=0}^{N-1}\frac{1}{\sqrt{h}}|G_{\gamma}^{h/2}*(\chi^h(\cdot -nh) - \chi^h(\cdot -(n+1)h))|_{\mathbb{A}}^2.
\end{aligned}\nonumber
\end{align}
Recalling that $N = \alpha /\sqrt{h}$ we can rewrite the right hand side as

\begin{equation}
\frac{\alpha^2}{N}\sum_{n=0}^{N-1}\frac{1}{h\sqrt{h}}|G_{\gamma}^{h/2}*(\chi^h(\cdot -nh) - \chi^h(\cdot - (n+1)h))|^2_{\mathbb{A}}.
\end{equation}
This is an average of time shifts of $\alpha^2\frac{1}{h\sqrt{h}}|G_{\gamma}^{h/2}*(\chi^h - \chi^h(\cdot - h))|^2_{\mathbb{A}}$. Since $Nh = o(1)$ all these time shifts are small, thus the average has the same distributional limit as $\alpha^2\frac{1}{h\sqrt{h}}|G_{\gamma}^{h/2}*(\chi^h - \chi^h(\cdot - h))|^2_{\mathbb{A}}$. This proves (\ref{subst2}). The argument for (\ref{subst3}) is similar.

\subsection*{Existence of the normal velocities}

We now prove the existence of the normal velocities. Fix $1 \le i \le N$ and observe that for $w \in \{\gamma, \beta\}$ we have

\begin{equation}\label{ineqTrick}
\begin{split}
|\chi_i - \chi_i(-\tau)| \le &(\chi_i - \chi_i(-\tau)) G_w^h * (\chi_i - \chi_i(-\tau)) + |\chi_i - G_{w}^h * \chi_i| 
\\ &+ |\chi_i(-\tau) - G_w^h*\chi_i(-\tau)|,
\end{split}
\end{equation}
which follows simply by observing that $|\chi_i - \chi_i(\cdot-\tau)| = |\chi_i - \chi_i(\cdot-\tau)|^2 = (\chi_i - \chi_i(\cdot-\tau)G_w^h * (\chi_i - \chi_i(\cdot-\tau)) + (\chi_i - \chi_i(\cdot-\tau))(\chi - G_w^h * \chi) + (\chi_i(\cdot-\tau) - \chi_i)(\chi_i(\cdot-\tau) - G_w^h * \chi_i(\cdot-\tau))$. Using Jensen's inequality and the elementary identity (\ref{ineqAB}) in the Appendix we have

\begin{align}\label{puttinInt}
\begin{aligned}
|\chi_i - G_{w}^h * \chi_i| &\le \int G_{w}^h(z)|\chi_i(x) - \chi_i(x-z)| dz
\\ & = \int G_{w}^h(z)\chi_i(x)(1 - \chi_i(x-z)) dz + \int G_{w}^h(z)(1-\chi_i(x)) \chi_i(x-z)dz
\\ & = \sum_{k\neq i}\int G_{w}^h(z)\chi_i(x)\chi_k(x-z) dz + \sum_{k\neq i}\int G_{w}^h(z)\chi_k(x) \chi_i(x-z)| dz.
\end{aligned}
\end{align}
Now observe that by testing (\ref{provedConv}) with $G_{w}/K_{ij}$ (which is bounded, and thus admissible), we learn that 

\begin{equation}
\lim_{h\downarrow 0}\frac{1}{\sqrt{h}}\int G_{w}^h(z)\chi_i(x)\chi_k(x-z) dz = \int G_{w}(z)(\nu_{ik}(x,t) \cdot z)_{+}dz\Hm_{|\Sigma_{ik}(t)}(dx)dt.
\end{equation}
Thus, if we divide (\ref{puttinInt}) by $\sqrt{h}$ and let $h \downarrow 0$, using also (\ref{DissMeasBoundOneInt}) we obtain

\begin{equation}\label{derivativeBound}
\begin{split}
\alpha|\partial_t\chi_i| &\le \liminf_{h\downarrow 0}\frac{|\delta\chi_i|}{\sqrt{h}} 
\\ &\le \limsup_{h\downarrow 0}\frac{|\delta\chi_i|}{\sqrt{h}} 
\\ &\le C\alpha^2\omega + C\Hm_{\partial^* \Omega_i(t)}(dx)dt,
\end{split}
\end{equation}
where $C$ is a constant which depends on $\gamma, \beta, N$, the mobilities and the surface tensions. If we divide by $\alpha$ and then let $\alpha \to 0$ we learn that $|\partial_t\chi_i|$ is absolutely continuous with respect to $\Hm_{|\partial^*\Omega_i(t)}(dx)dt$. In particular, there exists $V_i \in L^1(\Hm_{|\partial^*\Omega_i(t)}(dx)dt)$ which is the normal velocity of $\chi_i$ in the sense that $\partial_t\chi_i = V_i|\nabla \chi_i|$ in the sense of distributions. The optimal integrability $V_i \in L^2(\Hm_{|\partial^*\Omega_i(t)}(dx)dt)$ will be shown in the second part of the proof. Let us record for later use that with a similar reasoning we actually obtain that $\limsup_{h} \frac{|\delta\chi_i|}{\sqrt{h}}$ is absolutely continuous with respect to $\Hm_{|\partial^*\Omega_i(t)}(dx)dt$.Thus in particular inequality (\ref{derivativeBound}) holds with $\omega$ replaced by its absolutely continuous part with respect to $\Hm_{|\partial^*\Omega_i(t)}(dx)dt$; calling this $\omega^{ac}_i$, it means

\begin{equation}\label{errorBoundsDeltas}
\limsup_{h\downarrow 0} \frac{|\delta\chi_i|}{\sqrt{h}} \le C\alpha^2\omega_i^{ac} + C\Hm_{|\partial^*\Omega_i(t)}(dx)dt.
\end{equation}

\subsection*{Sharp Bound}

Before entering into the proof of the sharp bound, we need to prove the following property. For any $i \neq j$ we have that, in a distributional sense, the following holds

\begin{equation}\label{PlusPlus}
\lim_{h \downarrow 0} \frac{1}{\sqrt{h}} \delta\chi_i^+ K_{ij}^h*\delta\chi_j^+ = 0 = \lim_{h \downarrow 0} \frac{1}{\sqrt{h}} \delta\chi_i^- K_{ij}^h*\delta\chi_j^-.
\end{equation}

We focus on the first limit, the second one being analogous. The first observation is that the limit

\begin{equation}\label{lambdaBoh}
\lambda := \lim_{h\downarrow 0} \frac{1}{\sqrt{h}}\delta\chi_i^+K_{ij}^h * \delta\chi_j^+
\end{equation}
is a nonnegative bounded measure, which is absolutely continous with respect to $\Hm_{|\Sigma_{ij}(t)}(dx)dt$. Indeed, spelling out the $z$ integral and using the fact that $\delta\chi_i^+ = \chi_i(1-\chi_i(-\tau))$ we obtain

\begin{align}
&\begin{aligned}
\frac{1}{\sqrt{h}}\delta\chi_i^+K_{ij}^h * \delta\chi_j^+ = \frac{1}{\sqrt{h}}\int K_{ij}^h(z)\delta\chi_i^+(x,t) * \delta\chi_j^+(x-z, t)dz
\end{aligned}
\\ &\begin{aligned}
\le \frac{1}{\sqrt{h}}\int K_{ij}^h(z)\chi_i(x,t)\chi_j(x-z, t)dz
\end{aligned}
\end{align}
which by (\ref{provedConv}) in Lemma \ref{compactnessLemma}, as $h\downarrow 0$, converges to 

\begin{equation}
\int K_{ij}(z) (\nu_{ij}(x,t)\cdot z)_{+} \Hm_{|\Sigma_{ij}(t)}(dx)dt
\end{equation}
which is absolutely continous with respect to $\Hm_{|\Sigma_{ij}(t)}(dx)dt$. 

Now, given $\nu_0 \in \mathbf{S}^{d-1}$ we claim that

\begin{align}
\begin{aligned}\label{crucialBoundTrick}
\lambda \le &\int_{\nu_0 \cdot z \le 0} K_{ij}(z)(\nu_{ij}\cdot z)_+\Hm_{|\Sigma_{ij}(t)}(dx)dt 
\\ &+ \int_{\nu_0 \cdot z \ge 0} K_{ij}(z)(\nu_{ij}\cdot z)_-\Hm_{|\Sigma_{ij}(t)}(dx)dt.
\end{aligned}
\end{align}
To see this, let us denote momentarily the right-hand side of  (\ref{lambdaBoh}) (disintegrated in the $z$-variable) as $\lambda_h :=  \chi_i(x,t)(1-\chi_i)(x, t-\tau)K_{ij}^h(z)\chi_i(x-z,t)(1-\chi_i)(x-z,t-\tau)$. Using the fact that $0 \le \chi_i, \chi_j \le 1$ and $\sum_l \chi_l = 1$ we obtain the following inequalities

\begin{align}
&\begin{aligned}\label{ineqLE}
\lambda_h\le \chi_i(x,t)K_{ij}^h(z)\chi_i(x-z,t).
\end{aligned}
\\ &\begin{aligned}\label{ineqGE}
\lambda_h \le &\chi_j(x,t-\tau)K_{ij}^h(z)\chi_i(x-z,t-\tau) 
\\ &+ C\sum_{k\neq i,j} K_{ij}^h(z)\left(|\delta\chi_k|(x,t) + |\delta\chi_k|(x-z,t)\right).
\end{aligned}
\end{align}
Here $C$ is a constant that does not depend on $h$. Using inequality (\ref{ineqLE}) on the domain $\{\nu_0 \cdot z \le 0\}$ and inequality (\ref{ineqGE}) on the domain $\{\nu_0 \cdot z \ge 0\}$ we obtain

\begin{align*}
\phantom{\lambda}
&\begin{aligned}
\mathllap{\lambda} \le & \limsup_{h\downarrow 0}\frac{1}{\sqrt{h}}\int_{\nu_0 \cdot z \le 0} \chi_i(x,t)K_{ij}^h(z)\chi_i(x-z,t)dz 
\\ & + \limsup_{h\downarrow 0}\frac{1}{\sqrt{h}}\int_{\nu_0 \cdot z \ge 0} \chi_j(x,t-\tau)K_{ij}^h(z)\chi_i(x-z,t-\tau)dz 
\\ & + C\sum_{k\neq i, j}\limsup_{h\downarrow 0}\left(\frac{1}{\sqrt{h}}\int K_{ij}^h(z)|\delta\chi_k|(x,t)dz + \frac{1}{\sqrt{h}}\int K_{ij}^h(z)|\delta\chi_k|(x-z,t)dz  \right)
\end{aligned}.
\end{align*}
Observe that for any $1\le k\le N$ we have

\begin{equation}
\limsup_{h\downarrow 0}\frac{1}{\sqrt{h}}\int K_{ij}^h(z)|\delta\chi_k|(x,t)dz = \frac{1}{\sqrt{h}}\int K_{ij}^h(z)|\delta\chi_k|(x-z,t)dz.
\end{equation}
This can be seen by showing that 

\begin{equation}\label{sameStuff}
\lim_{h\downarrow 0}\frac{1}{\sqrt{h}}\int K_{ij}^h(z)\left( |\delta\chi_k|(x,t)-|\delta\chi_k|(x-z,t) \right)dz
\end{equation}
which can be shown to be true by testing with an admissible test function, and putting the spatial shift $z$ on it. Thus recalling (\ref{provedConv}) and (\ref{errorBoundsDeltas}), we obtain that 

\begin{align}
\begin{aligned}
\lambda \le &\int_{\nu_0 \cdot z \le 0} K_{ij}(z)(\nu_{ij}\cdot z)_+\Hm_{|\Sigma_{ij}(t)}(dx)dt 
\\ &+ \int_{\nu_0 \cdot z \ge 0} K_{ij}(z)(\nu_{ij}\cdot z)_-\Hm_{|\Sigma_{ij}(t)}(dx)dt
\\ &+ C\sum_{k\neq i,j}\alpha^2\omega_k^{ac} + \Hm_{|\partial^*\Omega_k(t)}(dx)dt.
\end{aligned}
\end{align}
Since we already know that $\lambda$ is absolutely continous with respect to $\Hm_{|\Sigma_{ij}(t)}(dx)dt$, the same bound holds true if we replace the right hand side with its absolutely continuous part with respect to $\Hm_{|\Sigma_{ij}(t)}(dx)dt$. Observing that for $k \neq i, j$  by Lemma \ref{densInt} in the Appendix the measures $\Hm_{|\partial^*\Omega_k(t)}(dx)dt$ and $\Hm_{|\partial^*\Sigma_{ij}(t)}(dx)dt$ are mutually singular , this yields (\ref{crucialBoundTrick}). 

Writing $\lambda = \theta(x,t)\Hm_{|\Sigma_{ij}(t)}(dx)dt$ for some $L^1(\Hm_{|\Sigma_{ij}(t)}(dx)dt)$-function $\theta$ we obtain that inequality (\ref{crucialBoundTrick}) yields

\begin{align}\label{crucialBoundInt}
\begin{aligned}
\theta(x,t) \le \int_{\nu_0 \cdot z \le 0}&K_{ij}(z)(\nu_{ij}(x,t) \cdot z)_+dz
\\ &\int_{\nu_0 \cdot z \ge 0}K_{ij}(z)(\nu_{ij}(x,t) \cdot z)_-dz
\end{aligned}
\end{align}
for every $\nu_0 \in \mathbf{S}^{d-1}$ and $\Hm_{|\Sigma_{ij}(t)}(dx)dt$-a.e. $(x,t) \in [0,1)^d \times (0,T)$. By a separability argument, we see that the null set on which (\ref{crucialBoundInt}) does not hold can be chosen so that it is independent of the choice of $\nu_0$. If we select $\nu_0 = \nu_{ij}(x,t)$ this yields $\theta \le 0$ almost everywhere with respect to $\Hm_{|\Sigma_{ij}(t)}(dx)dt$. Since we already know that $\lambda$ is nonnegative this gives $\lambda = 0$.

Before getting the sharp bound, we also need to check that $V_{ij}$ is well defined, i.e.\ we need to prove that for any $i\neq j$ we have $V_i = -V_j$ a.e.\ with respect to $\Hm_{|\Sigma_{ij}(t)}(dx)dt$. To see this, we start by observing that if $\xi \in C^{\infty}_c([0,1)^d\times (0,T))$, thanks to the fact that $\sum_{k\neq i} \chi_k = 1 - \chi_i$, we get

\begin{align}
\begin{aligned}
\int \xi V_i\Hm_{|\partial^*\Omega_i(t)}(dx)dt &= -\int \partial_t\xi \chi_i dxdt
\\ &=\sum_{k\neq i} \int \partial_t\xi \chi_k dxdt
\\ &=-\sum_{k\neq i} \int \xi V_k\Hm_{|\partial^*\Omega_k(t)}(dx)dt.
\end{aligned} 
\end{align}
Choosing $\xi = f(t)g(x)$ for some $f \in C^{\infty}_c((0,T))$ and $g \in C^{\infty}([0,1)^d)$, by a separability argument, we obtain that for a.e.\ $t$ and every $g \in C^{\infty}([0,1)^d)$

\begin{align}\label{eqHolds}
\begin{aligned}
\int g V_i\Hm_{|\partial^*\Omega_i(t)}(dx) &= -\sum_{k\neq i} \int g V_k\Hm_{|\partial^*\Omega_k(t)}(dx).
\end{aligned} 
\end{align}

Pick $t$ such that (\ref{eqHolds}) holds. Let $g \in C^{\infty}([0,1)^d)$ and let $\rho_{m}$ be a partition of unity obtained by the construction of Section \ref{consPOU} applied to the function $\chi(\cdot, t)$ on the interface $\Sigma_{ij}(t)$. Then 

\begin{align}\label{eqPoU}
\begin{aligned}
\sum_{m \in \mathbf{N}}\int \rho_{m}g V_i\Hm_{|\partial^*\Omega_i(t)}(dx) &= -\sum_{m \in \mathbf{N}}\sum_{k\neq i} \int \rho_{m}g V_k\Hm_{\partial^*\Omega_k(t)}(dx).
\end{aligned} 
\end{align}
Passing to the limit $r\downarrow 0$ in (\ref{eqPoU}) we get by Lemma \ref{propertiesPOU} that

\begin{align}
\begin{aligned}
\int g V_i\Hm_{|\Sigma_{ij}(t)}(dx) &= -\int g V_j\Hm_{\Sigma_{ij}(t)}(dx).
\end{aligned}
\end{align}
Since this identity holds for any $g \in C^{\infty}([0,1)^d)$, a density argument gives $V_i(x, t) = -V_j(x,t)$ for $\Hm_{|\Sigma_{ij}(t)}$-a.e.\ $x$. In other words

\begin{equation}
\int |V_i(x,t) +V_j(x,t)|\Hm_{|\Sigma_{ij}(t)}(dx) = 0.
\end{equation}
Integrating in time yields that $V_i = -V_j$ a.e.\ with respect to $\Hm_{|\Sigma_{ij}(t)}(dx)dt$.

We now proceed with the derivation of the sharp lower bound. Define $c_{ij} := \int K_{ij}(z)dz$. Then we have

\begin{align}\label{identityInterface}
&\begin{aligned}
c_{ij}(|\delta\chi_i| + |\delta\chi_j|) =  c_{ij}(\delta\chi_i^+ + \delta\chi_j^- + \delta\chi_i^- + \delta\chi_j^+)
\end{aligned}
\\ &\begin{aligned}
=\frac{1}{2}&\left( \delta\chi_i^+ K_{ij}^h * (1 - \delta\chi_j^-) + (1-\delta\chi_j^-) K_{ij}^h * \delta\chi_i^+ + \delta\chi_j^- K_{ij}^h * (1 - \delta\chi_i^+) \right.
\\ & \left.+ (1-\delta\chi_i^+) K_{ij}^h * \delta\chi_j^- + \delta\chi_i^- K_{ij}^h * (1 - \delta\chi_j^+) + (1-\delta\chi_j^+) K_{ij}^h * \delta\chi_i^- \right.
\\ & \left. + \delta\chi_j^+ K_{ij}^h * (1 - \delta\chi_i^-) + (1-\delta\chi_i^-) K_{ij}^h * \delta\chi_j^+\right) + \left(\delta\chi_i^+ K_{ij}^h*\delta\chi_j^- \right.
\\ & \left. + \delta\chi_j^- K_{ij}^h*\delta\chi_i^+ + \delta\chi_i^- K_{ij}^h*\delta\chi_j^+ + \delta\chi_j^+ K_{ij}^h*\delta\chi_i^- \right)
\end{aligned}\nonumber
\end{align}
Now we rewrite the terms in the second parenthesis using $-ab = a_+b_- + a_-b_+ - a_+b_+ - a_-b_-$ and then adding and subtracting the contributions of the minority phases we obtain
\begin{align}\label{NewCrucBd}
\begin{aligned}
c_{ij}(|\delta\chi_i| + |\delta\chi_j|)  \le \frac{1}{2}&\left( \delta\chi_i^+ K_{ij}^h * (1 - \delta\chi_j^-) + (1-\delta\chi_j^-) K_{ij}^h * \delta\chi_i^+ + \delta\chi_j^- K_{ij}^h * (1 - \delta\chi_i^+) \right.
\\ & \left. + (1-\delta\chi_i^+) K_{ij}^h * \delta\chi_j^- + \delta\chi_i^- K_{ij}^h * (1 - \delta\chi_j^+) + (1-\delta\chi_j^+) K_{ij}^h * \delta\chi_i^- \right.
\\ & \left. + \delta\chi_j^+ K_{ij}^h * (1 - \delta\chi_i^-) + (1-\delta\chi_i^-) K_{ij}^h * \delta\chi_j^+\right) - \sum_{l,p}\delta\chi_l K_{lp}^h*\delta\chi_p
\\ & +\delta\chi_i^+ K_{ij}^h*\delta\chi_j^+ + \delta\chi_i^- K_{ij}^h*\delta\chi_j^- + \delta\chi_j^+ K_{ij}^h*\delta\chi_i^+
\\ &+\delta\chi_j^- K_{ij}^h*\delta\chi_i^- + \sum_{\{l,p\} \neq \{i,j\}, \{l,p\}}\delta\chi_l K_{lp}^h*\delta\chi_p.
\end{aligned}
\end{align}

Now the main idea is to split the integral of $K_{ij}$ in the definition of $c_{ij}$ into two parts. More precisely, by the symmetry (\ref{symmetryHK}), for any $\nu_0 \in \mathbf{S}^{d-1}$ and any $V_0 > 0$ we have

\begin{equation}
c_{ij} = 2\int_{0 \le \nu_0 \cdot z \le \alpha V_0}K_{ij}(z)dz + 2\int_{\nu_0 \cdot z > \alpha V_0}K_{ij}(z)dz.
\end{equation}
Substituting into (\ref{NewCrucBd}) and dividing by $\sqrt{h}$ we obtain

\begin{align}\label{splitIdentity}
&\begin{aligned}
2\int_{0 \le \nu_0 \cdot z \le \alpha V_0} K_{ij}(z)dz\frac{(|\delta\chi_i| + |\delta\chi_j|)}{\sqrt{h}}
\end{aligned}
\\ &\begin{aligned}
= \frac{1}{2\sqrt{h}}&\bigg( \delta\chi_i^+ K_{ij}^h * (1 - \delta\chi_j^-) + (1-\delta\chi_j^-) K_{ij}^h * \delta\chi_i^+ + \delta\chi_j^- K_{ij}^h * (1 - \delta\chi_i^+)
\\ &+ (1-\delta\chi_i^+) K_{ij}^h * \delta\chi_j^- + \delta\chi_i^- K_{ij}^h * (1 - \delta\chi_j^+) + (1-\delta\chi_j^+) K_{ij}^h * \delta\chi_i^-
\\ &+ \delta\chi_j^+ K_{ij}^h * (1 - \delta\chi_i^-) + (1-\delta\chi_i^-) K_{ij}^h * \delta\chi_j^+
\\ &-4\int_{\nu_0 \cdot z > \alpha V_0}K_{ij}(z)dz(|\delta\chi_i| + |\delta\chi_j|) 
\\ &- 2\sum_{lp}\delta\chi_l K_{lp}^h*\delta\chi_p + 2\delta\chi_i^+ K_{ij}^h*\delta\chi_j^+ + 2\delta\chi_i^- K_{ij}^h*\delta\chi_j^-
\\ &+ 2\delta\chi_j^+ K_{ij}^h*\delta\chi_i^+ + 2\delta\chi_j^- K_{ij}^h*\delta\chi_i^-
\\ &+ \sum_{(l,p) \neq (i,j), (l,p) \neq (j,i)}\delta\chi_l K_{lp}^h*\delta\chi_p \bigg).
\end{aligned}\nonumber
\end{align}

We will be interested in bounding the $\liminf$ of the left hand side. Observe that the distributional limit of the last five terms is non-positive. Indeed, the limit of first four terms vanish distributionally by property (\ref{PlusPlus}), while the last term is bounded from above by

\begin{equation}
 \sum_{(l,p) \neq (i,j), (l,p) \neq (j,i)}\delta\chi_l^+ K_{lp}^h*\delta\chi_p^+ + \delta\chi_l^- K_{lp}^h*\delta\chi_p^-,\nonumber
\end{equation}
which vanish distributionally by property (\ref{PlusPlus}). We thus obtain that the $\liminf$ of the left hand side of (\ref{splitIdentity}) is bounded from above by

\begin{align}
\begin{aligned}
\liminf_{h\downarrow 0}\frac{1}{2\sqrt{h}}&\bigg( \delta\chi_i^+ K_{ij}^h * (1 - \delta\chi_j^-) + (1-\delta\chi_j^-) K_{ij}^h * \delta\chi_i^+ + \delta\chi_j^- K_{ij}^h * (1 - \delta\chi_i^+)
\\ &+ (1-\delta\chi_i^+) K_{ij}^h * \delta\chi_j^- + \delta\chi_i^- K_{ij}^h * (1 - \delta\chi_j^+) + (1-\delta\chi_j^+) K_{ij}^h * \delta\chi_i^-
\\ &+ \delta\chi_j^+ K_{ij}^h * (1 - \delta\chi_i^-) + (1-\delta\chi_i^-) K_{ij}^h * \delta\chi_j^+
\\ &- 4\int_{\nu_0 \cdot z > \alpha V_0}K_{ij}(z)dz(|\delta\chi_i| + |\delta\chi_j|) - 2\sum_{lp}\delta\chi_l K_{lp}^h*\delta\chi_p\bigg).
\end{aligned}
\end{align}

For the last term we use the sharp bound (\ref{DissMeasBound}), relating this term to our dissipation measure $\omega$. We would like to get a good bound for the other terms. This cannot be done naively as before, since we want the bound to be sharp. We claim that

\begin{align}\label{sharpBoundWanted}
&\begin{aligned}
\limsup_{h\downarrow 0}\frac{1}{\sqrt{h}}&\left( \delta\chi_i^+ K_{ij}^h * (1 - \delta\chi_j^-) + (1-\delta\chi_j^-) K_{ij}^h * \delta\chi_i^+ \delta\chi_j^- K_{ij}^h * (1 - \delta\chi_i^+) \right.
\\ & \left. + (1-\delta\chi_i^+) K_{ij}^h * \delta\chi_j^- + \delta\chi_i^- K_{ij}^h * (1 - \delta\chi_j^+) + (1-\delta\chi_j^+) K_{ij}^h * \delta\chi_i^- \right.
\\ & \left. + \delta\chi_j^+ K_{ij}^h * (1 - \delta\chi_i^-) + (1-\delta\chi_i^-) K_{ij}^h * \delta\chi_j^+ \right.
\\ &\left.- 4\int_{\nu_0 \cdot z > \alpha V_0}K_{ij}(z)dz(|\delta\chi_i| + |\delta\chi_j|) \right)
\end{aligned}
\\ &\begin{aligned} \le 8\int_{0 \le \nu_0 \cdot z \le \alpha V_0} K_{ij}(z)|\nu_{ij}(x)\cdot z|&dz \Hm_{|\Sigma_{ij}(t)}(dx)dt \\ & + C\sum_{k\neq i,j}(\alpha^2 \omega^{ac}_k + \Hm_{|\partial^*\Omega_k(t)}(dx)dt.
\end{aligned}\nonumber
\end{align}
Here $C$ is a constant that depends on $\gamma, \beta, \mathbb{A}, \mathbb{B}$, but not on $h$. Assume for the moment that (\ref{sharpBoundWanted}) is true and let us conclude the argument in this case. Using (\ref{sharpBoundWanted}) and (\ref{DissMeasBound}) we obtain

\begin{align}\label{distributionalBound}
&\begin{aligned}
2\liminf_{h\downarrow 0}\int_{0 \le \nu_0 \cdot z \le \alpha V_0} K_{ij}(z)dz\frac{(|\delta\chi_i| + |\delta\chi_j|)}{\sqrt{h}}
\end{aligned}
\\ &\begin{aligned} \le \alpha^2\omega + 4\int_{0 \le \nu_0 \cdot z \le \alpha V_0} K_{ij}(z)|\nu_{ij}(x)\cdot z|&dz \Hm_{|\Sigma_{ij}(t)}(dx)dt 
\\ & + C\sum_{k \neq i,j}(\alpha^2 \omega^{ac}_k + \Hm_{|\partial^*\Omega_k(t)}(dx)dt
\end{aligned}\nonumber
\end{align}
in the sense of distributions on $[0,1)^d\times (0,T)$. Observe also that the left hand side of (\ref{distributionalBound}) is an upper bound for $\int_{0 \le \nu_0 \cdot z \le \alpha V_0}K_{ij}(z)dz(|\partial_t\chi_i| + |\partial_t\chi_j|)$, thus the inequality still holds true if the left hand side is replaced by this term. Remember that $\omega_k^{ac}$ is absolutely continuous with respect to $\Hm_{|\partial^*\Omega_k(t)}(dx)dt$, thus there exist functions $W_k \in L^1(\Hm_{|\partial^*\Omega_k(t)}(dx)dt)$ such that $\omega_k^{ac} = W_k(x,t)\Hm_{|\partial^*\Omega_k(t)}(dx)dt$. We now disintegrate the measure $\omega$, i.e.\ we find a Borel family $\omega_{t}, t \in (0,T)$, of positive measures on $[0,1)^d$ such that $\omega = \omega_t \otimes dt$. Having said this, it is not hard to see that (\ref{distributionalBound}) holds in a disintegrated version, i.e.\ we have for Lebesgue a.e.\ $t \in (0,T)$

\begin{align}\label{distributionalDisintegratedBound}
&\begin{aligned}
2\int_{0 \le \nu_0 \cdot z \le \alpha V_0}K_{ij}(z)dz(|V_i(x,t)|\Hm_{|\partial^*\Omega_i(t)}(dx)+|V_j(x,t)|\Hm_{|\partial^*\Omega_j(t)}(dx)) 
\end{aligned}
\\ &\begin{aligned}\le \alpha^2\omega_t + 4\int_{0 \le \nu_0 \cdot z \le \alpha V_0} K_{ij}(z)|\nu_{ij}(x)\cdot z|&dz \Hm_{|\Sigma_{ij}(t)}(dx) 
\\ &+ C\sum_{k \neq i,j}(\alpha^2 W_k(x,t) +  1) \Hm_{|\partial^*\Omega_k(t)}(dx).
\end{aligned}\nonumber
\end{align}
Here $\nu_0 \in \mathbf{S}^{d-1}$ and $V_0 \in (0, \infty)$ are arbitrary: indeed even if the set of points in time for which (\ref{distributionalDisintegratedBound}) holds is a priori dependent on $\nu_0$ and $V_0$, a standard separability argument allows us to conclude that we can get rid of this dependence. 

Fix a point $t$ in time such that (\ref{distributionalDisintegratedBound}) holds. In what follows, we drop the time variable $t$ which is fixed, so for example $V_i(x) = V_i(x,t)$, $\Sigma_{ij} = \Sigma_{ij}(t)$ and so on. Fix $\xi \in C([0,1)^d)$, observe that by definition of $V_{ij}$ and by using the fact that $\Sigma_{ij} \subset \partial^*\Omega_i \cap \partial^*\Omega_j$ we have

\begin{align}
&\begin{aligned}
4\alpha\int_{0 \le \nu_0 \cdot z \le \alpha V_0}K_{ij}(z)dz \int_{[0,1)^d} \xi(x)|V_{ij}(x)|\Hm_{|\Sigma_{ij}}(dx)
\end{aligned}
\\ &\begin{aligned} 
\le\alpha^2\int_{[0,1)^d}\xi(x)\omega_t(dx)&+4\int_{[0,1)^d}\int_{0 \le \nu_0 \cdot z \le \alpha V_0} K_{ij}(z)|\nu_{ij}(x)\cdot z|dz\xi(x) \Hm_{|\Sigma_{ij}(t)}(dx) 
\\ & + C\sum_{k \neq i,j}\int_{[0,1)^d}\xi(x)(\alpha^2 W_k(x,t) +  1) \Hm_{|\partial^*\Omega_k(t)}(dx).
\end{aligned}\nonumber
\end{align}

Let us relabel $\nu_0$, $V_0$ and $\xi$ to make clear that they may depend on the pair $i,j$. Thus $\nu_0^{ij} \in \mathbf{S}^{d-1}$, $V_0^{ij} \in (0, \infty)$ and $\xi_{ij} \in C([0,1)^d)$ are arbitrary, and it holds

\begin{align}
&\begin{aligned}
2\alpha\int_{0 \le \nu_0^{ij} \cdot z \le \alpha V^{ij}_0}K_{ij}(z)dz \int_{[0,1)^d} \xi_{ij}(x)|V_{ij}(x)|\Hm_{|\Sigma_{ij}}(dx)
\end{aligned}
\\ &\begin{aligned} 
\le\alpha^2\int_{[0,1)^d}\xi_{ij}(x)\omega_t(dx)&+4\int_{[0,1)^d}\int_{0 \le \nu_0^{ij} \cdot z \le \alpha V_0^{ij}} K_{ij}(z)|\nu_{ij}(x)\cdot z|dz\xi_{ij}(x) \Hm_{|\Sigma_{ij}(t)}(dx) 
\\ & + C\sum_{k \neq i,j}\int_{[0,1)^d}\xi_{ij}(x)(\alpha^2 W_k(x,t) +  1) \Hm_{|\partial^*\Omega_k(t)}(dx).
\end{aligned}\nonumber
\end{align}

Let $\{\rho_{m}\}$ be a partition of unity obtained using the construction of Section \ref{consPOU} applied to the function $\chi(\cdot, t)$ on the inferface $\Sigma_{ij}(t)$. Use the above inequality with $\xi_{ij}$ replaced by $\rho_{m}\xi_{ij}$ and sum over $m$ and $i, j$ to get

\begin{align}
&\begin{aligned}
\sum_{i<j}\sum_{m\in \mathbf{N}} \mathbf{LH}_m^{ij} \le \sum_{i<j}\sum_{m \in \mathbf{N}} (\mathbf{I}_m^{ij}+\mathbf{II}_m^{ij}+\mathbf{III}_m^{ij})
\end{aligned}
\end{align}
where we have set

\begin{align}
\phantom{\mathbf{III}^{ij}_m}
&\begin{aligned}
\mathllap{\mathbf{LH}_m^{ij}} = 2\alpha\int_{0 \le \nu_0^{ij} \cdot z \le \alpha V^{ij}_0}K_{ij}(z)dz \int_{[0,1)^d} \rho_{mij}(x)\xi_{ij}(x)|V_{ij}(x)|\Hm_{|\Sigma_{ij}}(dx),
\end{aligned}\nonumber
\\ &\begin{aligned}
\mathllap{\mathbf{I}^{ij}_m} = \alpha^2\int_{[0,1)^d}\xi_{ij}(x)\rho_{mij}(x)\omega_t(dx),
\end{aligned}\nonumber
\\ &\begin{aligned}
\mathllap{\mathbf{II}^{ij}_m} = 4\int_{[0,1)^d}\int_{0 \le \nu_0^{ij} \cdot z \le \alpha V_0^{ij}} K_{ij}(z)|\nu_{ij}(x)\cdot z|dz\rho_{mij}(x)\xi_{ij}(x) \Hm_{|\Sigma_{ij}(t)}(dx),
\end{aligned}\nonumber
\\ &\begin{aligned}
\mathllap{\mathbf{III}^{ij}_m} = C\sum_{k \neq i,j}\int_{[0,1)^d}\rho_{mij}(x)\xi_{ij}(x)(\alpha^2 W_k(x,t) +  1) \Hm_{|\partial^*\Omega_k(t)}(dx).
\end{aligned}\nonumber
\end{align}
Observe that 

\begin{equation}
\sum_{i<j}\sum_{m \in \mathbf{N}} \mathbf{I}_m^{ij} \le \sum_{i<j} \alpha^2\int_{[0,1)^d}\xi_{ij}(x)\omega_t(dx)
\end{equation}
because $\rho_{m}$ is a partition of unity. Moreover by Lemma \ref{propertiesPOU} we get

\begin{align}
\phantom{\lim_{r\downarrow 0} \sum_{i<j}\sum_{m \in \mathbf{N}}\mathbf{III}_m^{ij}}
&\begin{aligned}
\mathllap{\lim_{r\downarrow 0}\sum_{i<j}\sum_{m \in \mathbf{N}}\mathbf{LH}_m^{ij}} =2\alpha\int_{0 \le \nu_0^{ij} \cdot z \le \alpha V^{ij}_0}K_{ij}(z)dz \int_{[0,1)^d}\xi_{ij}(x)|V_{ij}(x)|\Hm_{|\Sigma_{ij}}(dx).
\end{aligned}\nonumber
\\ &\begin{aligned}
\mathllap{\lim_{r\downarrow 0} \sum_{i<j}\sum_{m \in \mathbf{N}} \mathbf{II}_m^{ij}} =\sum_{i<j} 4\int_{[0,1)^d}\int_{0 \le \nu_0^{ij} \cdot z \le \alpha V_0^{ij}} K_{ij}(z)|\nu_{ij}(x)\cdot z|dz\xi_{ij}(x) \Hm_{|\Sigma_{ij}(t)}(dx). 
\end{aligned}\nonumber
\\ &\begin{aligned}
\mathllap{\lim_{r\downarrow 0} \sum_{i<j}\sum_{m \in \mathbf{N}} \mathbf{III}_m^{ij}} = 0.
\end{aligned}\nonumber
\end{align}

Putting things together we obtain that for any $\nu_0^{ij} \in \mathbf{S}^{d-1}$, any $V_0^{ij} \in (0, \infty)$, and any $\xi \in C([0,1)^d)$

\begin{align}\label{inequalityToApr}
&\begin{aligned}2\alpha\sum_{i<j}\int_{0 \le \nu_0^{ij} \cdot z \le \alpha V^{ij}_0}K_{ij}(z)dz \int_{[0,1)^d}\xi_{ij}(x)|V_{ij}(x)|\Hm_{|\Sigma_{ij}}(dx) 
\end{aligned}
\\ &\begin{aligned} \le &\sum_{i<j}\alpha^2\int_{[0,1)^d}\xi_{ij}(x)\omega_t(dx) 
\\ &+  4\sum_{i<j}\int_{[0,1)^d}\int_{0 \le \nu_0^{ij} \cdot z \le \alpha V_0^{ij}} K_{ij}(z)|\nu_{ij}(x)\cdot z|dz\xi_{ij}(x) \Hm_{|\Sigma_{ij}(t)}(dx).
\end{aligned}\nonumber
\end{align}

We now claim that by approximation the above inequality is valid for any simple function $\xi_{ij} \ge 0$. To see this, it is clear that we can concentrate on $\xi_{ij} = w_{ij}\mathbf{1}_{B_{ij}}$, where $B_{ij} \subset [0,1)^d$ are Borel and $w_{ij} \ge 0$. Observe that by the dominated convergence theorem, the family

\begin{equation}
\mathcal{F} := \left\{ B = \prod_{i<j} B_{ij}:\ B_{ij} \in \mathcal{B}([0,1)^d)\ \text{s.t.}\ \forall w_{ij} \ge 0\ \text{ (\ref{inequalityToApr}) holds with}\ \xi_{ij} = w_{ij} \mathbf{1}_{B_{ij}}\right\}
\end{equation}
is a monotone class. Thus by the monotone class theorem we just need to show that it contains all the products of open sets. But this is easy because given $B_{ij} \subset [0,1)^d$ open sets, we can always find sequences $\eta_{k}^{ij}$ of continuous functions with compact support such that $0 \le \eta_k^{ij} \le \mathbf{1}_{B_{ij}}$ and such that $\eta_k^{ij} \to \mathbf{1}_{B_{ij}}$, thus the claim follows by the monotone convergence theorem. 

With this in place one can use an approximation argument to replace the vector $\nu_0^{ij}$ with the $\Hm$-measurable vector valued function $\nu_{ij}$ obtaining the following inequality:

\begin{align}\label{ineqBefLim}
&\begin{aligned}2\alpha\sum_{i<j}\int_{[0,1)^d}\int_{0 < \nu_{ij}(x) \cdot z < \alpha V^{ij}_0}K_{ij}(z)dz \xi_{ij}(x)|V_{ij}(x)|\Hm_{|\Sigma_{ij}}(dx) \le
\end{aligned}
\\ &\begin{aligned} \le &\sum_{i<j} \alpha^2\int_{[0,1)^d}\eta_{ij}(x)\omega_t(dx)
\\ & +  4\sum_{i<j}\int_{[0,1)^d}\int_{0 \le \nu_{ij} \cdot z \le \alpha V_0^{ij}} K_{ij}(z)|\nu_{ij}(x)\cdot z|dz\xi_{ij}(x) \Hm_{|\Sigma_{ij}}(dx).
\end{aligned}\nonumber
\end{align}

Now divide by $\alpha^2$ and send $\alpha$ to zero. Record the following limits, which can be computed spelling out the definition of $K_{ij}$, and recalling the symmetry property (\ref{symmetryHK}) and the factorization property (\ref{factpropHK}) for the heat kernel

\begin{align}
&\begin{aligned}
\lim_{\alpha \downarrow 0} \frac{1}{\alpha}\int_{0 < \nu_{ij}(x) \cdot z < \alpha V_0^{ij}} K_{ij}(z) dz = \frac{V_0^{ij}}{\mu_{ij}}.
\end{aligned}
\\ &\begin{aligned}
\lim_{\alpha \downarrow 0} \frac{1}{\alpha^2}\int_{0 \le \nu_{ij}(x) \cdot z \le \alpha V_0^{ij}} K_{ij}(z)|\nu_{ij}(x)\cdot z| dz = \frac{(V_0^{ij})^2}{4\mu_{ij}}. 
\end{aligned}\nonumber
\end{align}
Then if we insert back into (\ref{ineqBefLim}) we obtain

\begin{align}\label{neededInCor}
&\sum_{i<j}\frac{2}{\mu_{ij}}\int_{[0,1)^d}V_0^{ij} \xi_{ij}(x)|V_{ij}(x)|\Hm_{|\Sigma_{ij}}(dx) 
\\ & \le \int_{[0,1)^d}\omega_t(dx) + \sum_{i<j}\int_{[0,1)^d}\frac{(V_0^{ij})^2}{\mu_{ij}}\xi_{ij}(x) \Hm_{|\Sigma_{ij}(t)}(dx).
\end{align}

Now given $M > 0$ take sequences of simple functions

\begin{equation}
s_m^{ij} = \sum_{k=1}^{p_m} w_k^{ij}\mathbf{1}_{B^k_{ij}}
\end{equation}
such that $s_m^{ij} \to |V_{ij}|\mathbf{1}_{\{|V_{ij}| \le M\}}$ as $m \to +\infty$ monotonically almost everywhere with respect to $\Hm_{|\Sigma_{ij}(t)}$. We are assuming that $\{B_{ij}^k\}_{k=1, . . . , p_m}$ are disjoint and $\Hm$-measurable,  with the property that $B_{ij}^{k_1} \cap B_{lr}^{k_2} = \emptyset$ if $\{i,j\} \neq \{l,r\}$. Choosing $V_0^{ij} = w_k^{ij}$, $\xi_{ij} = \mathbf{1}_{B^k_{ij}}$ in (\ref{neededInCor}) and summing over $k$ we obtain

\begin{align}
&\sum_{i<j}\frac{2}{\mu_{ij}}\int_{[0,1)^d}s_m^{ij}|V_{ij}(x)|\Hm_{|\Sigma_{ij}}(dx) 
\\ & \le \int_{[0,1)^d}\omega_t(dx) + \sum_{i<j}\int_{[0,1)^d}\frac{(s_m^{ij})^2}{\mu_{ij}} \Hm_{|\Sigma_{ij}(t)}(dx).\nonumber
\end{align}
Taking the limit $m\to +\infty$, using the monotone convergence theorem we obtain

\begin{align}
&\sum_{i<j}\frac{2}{\mu_{ij}}\int_{[0,1)^d}|V_{ij}(x)|^2\mathbf{1}_{\{|V_{ij}|\le M\}}\Hm_{|\Sigma_{ij}}(dx) 
\\ & \le \int_{[0,1)^d}\omega_t(dx) + \sum_{i<j}\int_{[0,1)^d}\frac{|V_{ij}(x)|^2}{\mu_{ij}}\mathbf{1}_{\{|V_{ij}|\le M\}}\Hm_{|\Sigma_{ij}}(dx)\nonumber
\end{align}
or, in other words,

\begin{align}
\sum_{i<j}\frac{1}{\mu_{ij}}\int_{[0,1)^d}|V_{ij}(x)|^2\mathbf{1}_{\{|V_{ij}|\le M\}}\Hm_{|\Sigma_{ij}}(dx) \le \int_{[0,1)^d}\omega_t(dx).\nonumber
\end{align}
Recall that $\mu_{ij} = \mu_{ji}$, thus the inequality above may be rewritten as

\begin{align}
\sum_{i,j}\frac{1}{2\mu_{ij}}\int_{[0,1)^d}|V_{ij}(x)|^2\mathbf{1}_{\{|V_{ij}|\le M\}}\Hm_{|\Sigma_{ij}}(dx) \le \int_{[0,1)^d}\omega_t(dx).\nonumber
\end{align}

If we now integrate in time we learn by the monotone convergence theorem that $V_{ij} \in L^2(\Hm_{|\Sigma_{ij}(t)}(dx)dt)$ and that the sharp bound (\ref{sharpVelocityBound}) is satisfied.

\subsection*{Proof of (\ref{sharpBoundWanted})}
To prove (\ref{sharpBoundWanted}) we proceed in several steps. 

First of all, we claim that the first eight terms may be substituted by

\begin{align}\label{halfSpaceInt}
2\int_{\nu_0 \cdot z \ge 0}K_{ij}^h(z)&\left( |\delta\chi_i^+ - \delta\chi_j^-(-z)| + |\delta\chi_i^+(-z) - \delta\chi_j^-|\right.
\\ &\left.|\delta\chi_i^- - \delta\chi_j^+(-z)| + |\delta\chi_i^-(-z) - \delta\chi_j^+| \right)dz.\nonumber
\end{align}

To show this, observe that we may replace the implicit $z$-integrals in the convolution in the first eight terms by twice the integrals over the half space $\{\nu_0 \cdot z \ge 0\}$ instead of $\mathbf{R}^d$. This is clearly true once we observe that

\begin{align}\label{integralHalfSpace}
&\begin{aligned}
\lim_{h\downarrow 0}\frac{1}{\sqrt{h}}&\left(\delta\chi_{i}^+ \int_{\nu_0\cdot z \ge 0} K_{ij}^h(z)(1-\delta\chi_j^-(\cdot - z)) dz \right.
\\ &\left. +(1- \delta\chi_{j}^-) \int_{\nu_0\cdot z \ge 0} K_{ij}^h(z)\delta\chi_i^+(\cdot - z) dz\right)
\end{aligned}
\\&\begin{aligned}
\mathllap{=}\lim_{h\downarrow 0}\frac{1}{\sqrt{h}}&\left(\delta\chi_{i}^+ \int_{\nu_0\cdot z \le 0} K_{ij}^h(z)(1-\delta\chi_j^-(\cdot - z)) dz \right.
\\ &\left.+(1- \delta\chi_{j}^-) \int_{\nu_0\cdot z \le 0} K_{ij}^h(z)\delta\chi_i^+(\cdot - z) dz\right)
\end{aligned}\nonumber
\end{align}
and that similar identities hold exchanging the roles  of $i,j$ and $+, -$ respectively. That (\ref{integralHalfSpace}) holds is not difficult to show. Indeed taking into account the fact that the kernel is even, the argument of the second limit is just a spatial shift of $z$ of the first one. The spatial shift may be put onto the test function, and thanks to the scaling of the kernel one can get the claim. We may thus substitute the first eight terms of the left hand side of (\ref{sharpBoundWanted}) with twice the same terms with the integration with respect to $z$ on the half space $\{\nu_0 \cdot z \ge 0\}$. If we rely again on the fact that $\delta\chi_i^+ \in \{0,1\}$, by identity (\ref{ineqAB}) in the Appendix we obtain (\ref{halfSpaceInt}), as claimed.

Now we need two inequalities for the integrand. First note that the integrand is a second-order finite difference, we claim that

\begin{align}\label{inequalitiesIntegrand}
&\begin{aligned}
|\delta\chi_i^+ - \delta\chi_j^-(\cdot-z)| + |\delta\chi_i^+(\cdot-z) - \delta\chi_j^-|+|\delta\chi_i^- - \delta\chi_j^+(\cdot-z)| + |\delta\chi_i^-(\cdot-z) - \delta\chi_j^+|
\end{aligned}
\\ &\begin{aligned}
\le \begin{cases}
|\delta\chi_i^+ - \delta\chi_i^+(\cdot-z)| + |\delta\chi_i^- - \delta\chi_i^-(\cdot-z)| + |\delta\chi_j^+ - \delta\chi_j^+(\cdot-z)| + |\delta\chi_j^- - \delta\chi_j^-(\cdot-z)| \\+ 4\sum_{k\neq i,j} (|\delta\chi_k| + \delta\chi_k(\cdot-z)|).
\\
\mathllap{ }
\\
|\delta\chi_i| + |\delta\chi_i(\cdot-z)| + |\delta\chi_j| + |\delta\chi_j(\cdot-z)|.
\end{cases}
\end{aligned}\nonumber
\end{align}

The second one follows from the triangle inequality. To show the first one, observe that

\begin{align}
\phantom{|\delta\chi_i^+ - \delta\chi_j^-(\cdot-z)|}
&\begin{aligned}
\mathllap{|\delta\chi_i^+ - \delta\chi_j^-(\cdot-z)|} = &(1-\delta\chi_i^+)\delta\chi_j^-(\cdot-z)  + \delta\chi_i^+(1-\delta\chi_j^-(\cdot-z))
\end{aligned}
\\ &\begin{aligned}
\mathllap{ }\le &(1-\delta\chi_i^+)\delta\chi_i^+(\cdot-z) + \sum_{k \neq i,j} |\delta\chi_k(\cdot-z)| + \delta\chi_j^+(1-\delta\chi_j^+(\cdot-z)) \\ &+ \sum_{k \neq i,j} |\delta\chi_k|
\end{aligned}\nonumber
\end{align}
and that similarly

\begin{align}
\phantom{|\delta\chi_i^+(\cdot-z) - \delta\chi_j^-|}
&\begin{aligned}
\mathllap{|\delta\chi_i^+(\cdot-z) - \delta\chi_j^-|} = &(1-\delta\chi_i(\cdot-z)^+)\delta\chi_j^- + \delta\chi_i(\cdot-z)^+(1-\delta\chi_j)
\end{aligned}
\\ &\begin{aligned}
\mathllap{ }\le &(1-\delta\chi_i(\cdot-z)^+)\delta\chi_i^- + \sum_{k \neq i,j} |\delta\chi_k| + \delta\chi_j(\cdot-z)^+(1-\delta\chi_j^-) \\ &+ \sum_{k \neq i,j} |\delta\chi_k(\cdot-z)|.
\end{aligned}\nonumber
\end{align}
Summing up the two inequalities we get

\begin{align}
&|\delta\chi_i^+ - \delta\chi_j^-(\cdot-z)| + |\delta\chi_i^+(\cdot-z) - \delta\chi_j^-| \le
\\ &\le |\delta\chi_i^+ - \delta\chi_i^+(\cdot-z)| + |\delta\chi_j^- - \delta\chi_j^-(\cdot-z)| + 2\sum_{k\neq i,j} |\delta\chi_k| + |\delta\chi_k(	\cdot-z)|.\nonumber
\end{align}
Similar bounds hold for the remaining terms in (\ref{inequalitiesIntegrand}). 

We now split the integral (\ref{halfSpaceInt}) into the domains of integration $\{0 \le \nu_0 \cdot z \le \alpha V_0\}$ and $\{ \nu_0 \cdot z > \alpha V_0\}$. On the first one we use the first inequality in (\ref{inequalitiesIntegrand}) for the integrand. Recalling identity (\ref{ineqAB}) and inequality (\ref{ineqLongAB}) in the Appendix we obtain, and using the fact that $\sum_k \chi_k = 1$

\begin{align}
&\begin{aligned}
2\int_{0 \le \nu_0 \cdot z \le \alpha V_0}K_{ij}^h(z)&\left( |\delta\chi_i^+ - \delta\chi_j^-(\cdot-z)| + |\delta\chi_i^+(\cdot-z) - \delta\chi_j^-|\right.
\\ &\left.+|\delta\chi_i^- - \delta\chi_j^+(\cdot-z)| + |\delta\chi_i^-(\cdot-z) - \delta\chi_j^+| \right)dz
\end{aligned}\nonumber
\\ &\begin{aligned}
\le 2\int_{0 \le \nu_0 \cdot z \le \alpha V_0}K_{ij}^h(z)&\left( |\chi_i - \chi_i(\cdot-z)| + |\chi_i(-\tau) - \chi_i(\cdot-\tau,\cdot-z)|\right.
\\ &\left.|\chi_j - \chi_j(\cdot-z)| + |\chi_j(\cdot-\tau) - \chi_j(\cdot-\tau,\cdot-z)| \right.
\\&\left. + 8\sum_{k\neq i,j} |\delta\chi_k| + |\delta\chi_k(\cdot-z)|\right)dz
\end{aligned}\nonumber
\\ & \begin{aligned}\label{tauTerms}
\le 2\int_{0 \le \nu_0 \cdot z \le \alpha V_0}K_{ij}^h(z)&\left( \chi_i\chi_j(\cdot-z) + \chi_i(\cdot-z)\chi_j + \sum_{k\neq i,j}\chi_i\chi_k(\cdot-z) + \chi_i(\cdot-z)\chi_k \right.
\\ &\left. \chi_i(\cdot-\tau)\chi_j(\cdot-\tau, \cdot-z) + \chi_i(\cdot-\tau,\cdot-z)\chi_j(\cdot-\tau) \right.
\\ & \left. + \sum_{k\neq i,j}\chi_i(\cdot-\tau)\chi_k(\cdot-\tau,\cdot -z) + \chi_i(\cdot-\tau, \cdot-z)\chi_k(\cdot-\tau) \right.
\\ & \left. \chi_j\chi_i(\cdot-z) + \chi_j(\cdot-z)\chi_i + \sum_{k\neq i,j}\chi_j\chi_k(\cdot-z) + \chi_j(\cdot-z)\chi_k \right.
\\ &\left. \chi_j(\cdot-\tau)\chi_i(\cdot-\tau, \cdot-z) + \chi_j(\cdot-\tau,\cdot-z)\chi_i(\cdot-\tau) \right.
\\ & \left. + \sum_{k\neq i,j}\chi_j(\cdot-\tau)\chi_k(\cdot-\tau, \cdot-z) + \chi_j(\cdot-\tau, \cdot-z)\chi_k(\cdot-\tau) \right.
\\&\left. + 8\sum_{k\neq i,j} |\delta\chi_k| + |\delta\chi_k(\cdot-z)|\right)dz.
\end{aligned}
\end{align}
On the set $\{\nu_0 \cdot z > \alpha V_0\}$ we use the second inequality in (\ref{inequalitiesIntegrand}), obtaining

\begin{align}\label{deltaTerms}
&\begin{aligned}
2\int_{ \nu_0 \cdot z > \alpha V_0} K_{ij}^h(z)&\left( |\delta\chi_i^+ - \delta\chi_j^-(\cdot-z)| + |\delta\chi_i^+(\cdot-z) - \delta\chi_j^-|\right.
\\ &\left.+|\delta\chi_i^- - \delta\chi_j^+(\cdot-z)| + |\delta\chi_i^-(\cdot-z) - \delta\chi_j^+| \right)dz
\end{aligned}
\\ &\begin{aligned}
\le 2\int_{ \nu_0 \cdot z > \alpha V_0} K_{ij}^h(z)(|\delta\chi_i| + |\delta\chi_i(\cdot-z)| + |\delta\chi_j| + |\delta\chi_j(\cdot-z)|)dz.
\end{aligned}\nonumber
\end{align} 
We now observe that for any $1 \le k \le N$ we have, as we already observed in (\ref{sameStuff})

\begin{equation}\label{shiftZ}
\lim_{h\downarrow 0} \frac{1}{\sqrt{h}}\int_{0 \le \nu_0 \cdot z \le \alpha V_0} K_{ij}^h(z)(|\delta\chi_k(\cdot-z)| - |\delta\chi_k|)dz = 0,
\end{equation}
thus in particular

\begin{align}\label{limsupShiftZ}
&\limsup_{h \downarrow 0} \frac{1}{\sqrt{h}}\int_{0 \le \nu_0 \cdot z \le \alpha V_0} K_{ij}^h(z)|\delta\chi_k(\cdot-z)|dz 
\\ &  = \limsup_{h\downarrow 0}\frac{1}{\sqrt{h}}\int_{0 \le \nu_0 \cdot z \le \alpha V_0} K_{ij}^h(z)|\delta\chi_k|dz.\nonumber
\end{align}
By putting the time shift $\tau$ on the test function if is easy to check that the distributional limit of the terms of (\ref{tauTerms}) which involve the shift $\tau$ have the same limit as the corresponding terms without the time shift. Thus recalling (\ref{errorBoundsDeltas}) and relying on (\ref{shiftZ}) and (\ref{provedConv}) we obtain that inserting (\ref{tauTerms}) and (\ref{deltaTerms}) into (\ref{halfSpaceInt}), the left hand side of (\ref{sharpBoundWanted}) is bounded by

\begin{align}
\begin{aligned}
&8\int_{0 \le \nu_0 \cdot z \le \alpha V_0}K_{ij}(z)((\nu_{ij}(x,t)\cdot z)_++(\nu_{ij}(x,t)\cdot z)_-)\Hm_{|\Sigma_{ij}(t)}(dx)dt 
\\ &+ C\sum_{k\neq i,j} \int_{0 \le \nu_0 \cdot z \le \alpha V_0}K_{ik}(z)((\nu_{ik}(x,t)\cdot z)_++(\nu_{ik}(x,t)\cdot z)_-)\Hm_{|\Sigma_{ik}(t)}(dx)dt
\\ &+C\sum_{k\neq i,j}(\alpha^2\omega_k^{ac} + \Hm_{|\partial^*\Omega_k(t)}(dx)dt),\nonumber
\end{aligned}
\end{align}
which clearly gives the claim once we realize that

\begin{align}
\begin{aligned}
&\int_{0 \le \nu_0 \cdot z \le \alpha V_0}K_{ik}(z)((\nu_{ik}(x)\cdot z)_++(\nu_{ik}(x)\cdot z)_-)\Hm_{|\Sigma_{ik}(t)}(dx)dt
\\ & \le 2\int_{\mathbf{R}^d}K_{ik}(z)|z|dz \le C.\nonumber
\end{aligned}
\end{align}
\end{proof}

\begin{proof}[Proof of Proposition \ref{meanCurv}]
The proof is along the same lines as Proposition 2 in \cite{Laux2019}, where the claim is analized in the case of two phases. For the convenience of the reader, we outline the strategy of the full proof, providing details only for the required changes. The proof is split into several steps. 
\\
\textsc{step 1}. The first observation is that for any $h>0$, any admissible $u \in \mathcal{M}$ and any smooth vector field $\xi$ we have the following lower bound for the metric slope, cf.\ (\ref{metricDerDef})

\begin{equation}
\frac{1}{2}|\partial E_h|(u) \ge \delta E_h(u)_\bullet\xi - \frac{1}{2}\left(\delta d_h(\cdot,u)_\bullet\xi\right)^2.\nonumber
\end{equation}
Here $\delta$ denotes the first variation, which is computed considering the curve $s \to u_s$ of configurations which solve the transport equations

\begin{equation}\label{transportEq}
\begin{cases}
\partial_su_i^s + \xi \cdot \nabla u_i^s = 0,
\\
u_i^s(\cdot, 0) = u_i(\cdot),
\end{cases}
\end{equation}
and by setting

\begin{equation}\label{defVar}
\delta E_h(u)_\bullet\xi := \frac{d}{ds}_{|{s=0}} E_h(u^s)\ \text{and}\  \delta d_h(\cdot, u)_\bullet\xi := \frac{d}{ds}_{|{s=0}} d(u, u^s).
\end{equation}
\\
\textsc{step 2}. The second observation is a representation formula for $\delta E_h(u)_\bullet\xi$. Namely

\begin{align}\label{repForE}
\begin{aligned}
\delta E_h(u)_\bullet\xi = \sum_{i,j}\frac{1}{\sqrt{h}}&\left(\int \nabla \cdot \xi u_iK_{ij}^h*u_jdx + \int \nabla \cdot \xi u_jK_{ij}^h*u_idxdt\right.
\\ &\left.+ \int [\xi, \nabla K_{ij}^h*](u_j)u_i dx\right).
\end{aligned}
\end{align}
Here $[\xi, \nabla K_{ij}^h*]$ denotes the commutator obtained taking the convolution with $\nabla K_{ij}^h$ and multiplying by $\xi$. To check this formula one starts by assuming $u$ to be smooth and then an approximation argument gives the result for a general $u \in \mathcal{M}$.
\\
\textsc{step 3}. Representation for $\delta d_h(\cdot, u)_\bullet\xi$. One checks that

\begin{align}
\begin{aligned}
&\frac{1}{2}\left(\delta d_h(\cdot, u)_\bullet\xi \right)^2
\\ &=\frac{\sqrt{h}}{2}\sum_{i,j}\left(\int u_i\xi\cdot \nabla^2K_{ij}^h *(\xi u_j)dx + \int u_j\xi\cdot \nabla^2K_{ij}^h *(\xi u_i)dx\right.
\\ &\left.+\int u_i \nabla\cdot \xi \nabla K_{ij}^h*(\xi u_j)dx + +\int u_j \nabla\cdot \xi \nabla K_{ij}^h*(\xi u_i)dx\right.
\\ &\left.-\int u_i \nabla\cdot \xi K_{ij}^h*(u_j\nabla\cdot \xi)dx-\int u_j \nabla\cdot \xi K_{ij}^h*(u_i\nabla\cdot \xi)dx\right.
\\ &\left.-\int \xi u_i \nabla K_{ij}^h*(u_j\nabla\cdot \xi)dx-\int \xi u_j \nabla K_{ij}^h*(u_i\nabla\cdot \xi)dx\right).\nonumber
\end{aligned}
\end{align}

Once again this formula can be easily checked when $u$ is smooth, an approximation argument then gives the extension to the case $u \in \mathcal{M}$.
\\
\textsc{step 4}. Passage to the limit in $\delta E_h$. We claim that

\begin{equation}\label{limitE}
\lim_{h\downarrow 0} \int_0^T \delta E_h(u^h(t))_\bullet\xi dt = \sum_{i,j} \sigma_{ij}\int \left(\nabla \cdot \xi - \nu_{ij} \cdot \nabla\xi \nu_{ij}\right) \Hm_{|\Sigma_{ij}(t)}(dx)dt.
\end{equation}

The proof is very similar to the two phases case, and relies on the weak convergence (\ref{provedConv}). Firstly, testing (\ref{provedConv}) with $\nabla \cdot \xi$ we get

\begin{align}
\begin{aligned}
&\lim_{h\downarrow 0} \sum_{i,j}\frac{1}{\sqrt{h}}\int \left( \nabla\cdot\xi u_i^hK_{ij}^h*u_j^h + \nabla\cdot\xi u_j^hK_{ij}^h*u_i^h \right)dxdt
\\ & =\sum_{i,j}2\sigma_{ij}\int \nabla\cdot\xi \Hm_{|\Sigma_{ij}(t)}(dx)dt.\nonumber
\end{aligned}
\end{align}
For the term involving the commutator, one checks that

\begin{align}
\lim_{h\downarrow 0} \left(\int [\xi,\nabla K_{ij}^h*](u_j^h)u_i^hdxdt - \int \nabla\xi z\cdot \nabla K_{ij}^h(z)u_j^h(x-z,t)u_i^h(x,t)dzdxdt \right) = 0.\nonumber
\end{align}
With this in place, we observe that 

\begin{align}
\begin{aligned}
&\int \nabla\xi z\cdot \nabla K_{ij}^h(z)u_j^h(x-z,t)u_i^h(x,t)dzdxdt
\\ & = \int \nabla \xi(x,t)z \cdot \nabla K_{ij}(z)(\nu_{ij}(x,t)\cdot z)_{+} \Hm_{|\Sigma_{ij}(t)}(dx)dt\nonumber
\end{aligned}
\end{align}
which can be seen by testing (\ref{provedConv}) with $\frac{\nabla\xi z \cdot \nabla K_{ij}(z)}{K_{ij}(z)}$ which is of polynomial growth in $z$. To conclude (\ref{limitE}) we just need to show that for any symmetric matrix $A \in \mathbf{R}^{d\times d}$ and any unit vector $\nu$ we have

\begin{equation}
\int Az \cdot \nabla K_{ij}(z)(\nu\cdot z)_{+}dz = -\sigma_{ij}\left( \operatorname{tr}A + \nu\cdot A\nu\right).\nonumber
\end{equation}
Using the definition of the kernel $K_{ij}$ it suffices to show that

\begin{equation}
\int Az \cdot \nabla G_{w}(z)(\nu\cdot z)_{+}dz = -\frac{\sqrt{w}}{\sqrt{\pi}}\left( \operatorname{tr}A + \nu\cdot A\nu\right)\ w \in \{\gamma, \beta\}.\nonumber
\end{equation}
\\
\textsc{step 5}. Passage to the limit in $\delta d_h(\cdot, u)\xi$. We claim that

\begin{align}\label{identityDH}
\begin{aligned}
\lim_{h\downarrow 0} \frac{1}{2}\left(\delta d_h(\cdot, u^h)_\bullet\xi\right) = \sum_{i,j}\frac{1}{2\mu_{ij}}\int (\xi \cdot \nu_{ij})^2 \Hm_{|\Sigma_{ij}(t)}(dx)dt.
\end{aligned}
\end{align}

To prove this, we observe that the terms which do not involve the Hessian $\nabla^2 K_{ij}^h$ are all $O(\sqrt{h})$. For example, to prove that 

\begin{equation}\label{provError}
\sqrt{h}\int u_i^h\nabla \cdot \xi \nabla K_{ij}^h *(\xi u_j^h)dxdt = O(\sqrt{h}),
\end{equation}
spell out the integral in the convolution, use the fact that $\nabla K_{ij}^h = \frac{1}{\sqrt{h}^{d+1}}\nabla K_{ij}(\frac{z}{\sqrt{h}})$, use the fact that $\nabla \xi(x,t)\xi(x-\sqrt{h}z,t)$ is bounded and test (\ref{provedConv}) with $\nabla K_{ij} / K_{ij}$. The other terms can be treated similarly. For the terms involving the Hessian of the kernel, we split the claim into

\begin{align}
&\begin{aligned}\label{mainPartTerm}
\lim_{h\downarrow 0}\sqrt{h}\int u_i^h(\xi \cdot \nabla^2K_{ij}^h *u_j)\xi dx dt = \frac{1}{2\mu_{ij}}\int (\xi \cdot \nu_{ij}(x,t))^2\Hm_{|\Sigma_{ij}(t)}(dx)dt,
\end{aligned}
\\ &\begin{aligned}\label{commutatorError}
\sqrt{h}\int u_i^h\xi\cdot [\xi, \nabla^2K_{ij}^h*](u_j^h)dxdt = O(\sqrt{h}).
\end{aligned}
\end{align}

The proof of (\ref{commutatorError}) is similar to the argument for (\ref{provError}). To prove identity (\ref{mainPartTerm}) observe that by spelling out the $z$-integral, a change of variable and by testing (\ref{provedConv}) with $\frac{\xi(x,t)\cdot\nabla^2 K_{ij}(z)\xi(x, t)}{K_{ij}(z)}$ we obtain

\begin{align}
\begin{aligned}
&\lim_{h\downarrow 0}\sqrt{h}\int u_i^h(\xi \cdot \nabla^2K_{ij}^h *u^j)\xi dx dt 
\\ &= \int \xi \cdot \nabla^2K_{ij}(z) \xi (\nu_{ij}(x,t)\cdot z)_+ \Hm_{|\Sigma_{ij}(t)}(dx)dt.\nonumber
\end{aligned}
\end{align}
Now identity (\ref{identityDH}) follows from the following formula: for any  two vectors $\xi \in \mathbf{R}^d$ and $\nu \in \mathbf{S}^{d-1}$ we have

\begin{equation}\label{secondFormula}
\int \xi \cdot \nabla^2K_{ij}(z) \xi (\nu\cdot z)_+ dz = \frac{1}{2\mu_{ij}}(\xi\cdot \nu)^2.
\end{equation}

To check (\ref{secondFormula}), by relying on the definition of the kernels, we just need to show that for $w \in \{\gamma, \beta\}$

\begin{equation}
\int \xi \cdot \nabla^2G_w(x) \xi (\nu\cdot z)_+ dz = \frac{1}{2\sqrt{\pi w}}(\xi\cdot \nu)^2.\nonumber
\end{equation}
Since the kernel is isotropic, we can reduce to the case $\xi = e_1$, thus we need to prove 

\begin{equation}
\int \partial^2_1G_w(x)(\nu\cdot z)_+ dz = \frac{1}{2\sqrt{\pi w}}\nu_1^2.\nonumber
\end{equation}
This can be done after two integration by parts and observing that 

\begin{equation}
\int_{\nu\cdot z = 0} G_w(z) dz = \frac{1}{2\sqrt{\pi w}}.\nonumber
\end{equation}

\textsc{conclusion}. By \textsc{step 1} we have

\begin{equation}
\frac{1}{2}\int_0^T|\partial E_h|^2(u^h)\ dt \ge \int_0^T \delta E_h(u^h)_\bullet\xi dt - \frac{1}{2}\int_0^T \left(\delta d_h(\cdot, u^h)_\bullet\xi\right)^2dt.\nonumber
\end{equation}
Taking the liminf on the left hand side, using \textsc{step 4} and \textsc{step 5} we get that for any smooth vector field $\xi$

\begin{align}
\begin{aligned}
\liminf_{h\downarrow 0} \frac{1}{2}\int_0^T|\partial E_h|^2(u^h) dt \ge &\sum_{i,j}\left[ \sigma_{ij}\int\left(\nabla \cdot \xi - \nu_{ij}\cdot \nabla\xi \nu_{ij}\right)\Hm_{|\Sigma_{ij}(t)}(dx)dt\right.
\\ &\left. - \frac{1}{2\mu_{ij}} \int (\xi\cdot \nu_{ij})^2\  \Hm_{|\Sigma_{ij}(t)}(dx)dt \right].\nonumber
\end{aligned}
\end{align}
Since the left hand side is bounded, the Riesz representation theorem for $L^2$ yields functions $H_{ij} \in L^2(\Hm_{|\Sigma_{ij}(t)}(dx)dt)$ such that 

\begin{equation}
\sum_{i,j}\sigma_{ij}\int\left(\nabla \cdot \xi - \nu_{ij}\cdot \nabla\xi \nu_{ij}\right)\ \Hm_{|\Sigma_{ij}(t)}(dx)dt =-\sum_{i,j}\sigma_{ij}\int H_{ij}\nu_{ij} \cdot \xi\ \Hm_{|\Sigma_{ij}(t)}(dx)dt\nonumber
\end{equation}
and such that for any $\xi \in L^2(\Hm_{|\bigcup_{i,j}\Sigma_{ij}(t)}(dx)dt)$

\begin{align}
\begin{aligned}
\liminf_{h\downarrow 0} \frac{1}{2}\int_0^T|\partial E_h|(u_h)\  dt \ge \sum_{i,j}\bigg( &-\sigma_{ij}\int H_{ij}\nu_{ij}\cdot\xi \Hm_{|\Sigma_{ij}(t)}(dx)dt
\\ &- \frac{1}{2\mu_{ij}} \int (\xi\cdot \nu_{ij})^2 \Hm_{|\Sigma_{ij}(t)}(dx)dt \bigg).\nonumber
\end{aligned}
\end{align}
Since the integration measures are mutually singular we can test with $\xi \in L^2(\Hm_{|\bigcup_{i,j}\Sigma_{ij}(t)}(dx)dt)$ such that $\xi_{|\Sigma_{ij}(t)} = -\mu_{ij}\sigma_{ij}H_{ij}\nu_{ij}$. This yields

\begin{align}
\begin{aligned}
\liminf_{h\downarrow 0} \frac{1}{2}\int_0^T|\partial E_h|^2(u_h)\ dt \ge &\sum_{i,j}\frac{\sigma_{ij}^2\mu_{ij}}{2}\int H_{ij}^2\ \Hm_{|\Sigma_{ij}(t)}(dx)dt.\nonumber
\end{aligned}
\end{align}
\end{proof}

\section{Appendix}\label{sec:appendix}

\subsection{Proof of Lemma \ref{propertiesPOU}}

Before giving the proof of this result, we need a simple technical lemma.

\begin{lemma}\label{densInt}
Fix $1 \le l \neq p \le N$. Then for any $1 \le i \neq j \le N$ such that $\{i,j\} \neq \{l,p\}$ the interfaces $\Sigma_{ij}$ and $\Sigma_{lp}$ are disjoint. In particular for $\Hm$-a.e. $x \in \Sigma_{lp}$ we have that

\begin{equation}\label{zerolimit}
\lim_{r \to 0}\Hm(\Sigma_{ij} \cap B(x,r)) = 0
\end{equation}
\end{lemma}
\begin{proof}
We first show that the interfaces $\Sigma_{ij}$ and $\Sigma_{lp}$ are disjoint. This follows immediately once we recall that every point in the reduced boundary of a set of finite perimeter has density $1/2$ (see \cite{Maggi2012}, Corollary 15.8). Assume for example that $i \neq l, p$. Thus if $y \in \Sigma_{lp}$ we have
\begin{equation}
\begin{split}
&1 \ge \limsup_r \frac{|(\Omega_l \cup \Omega_p \cup \Omega_i) \cap B(y, r)|}{\omega_dr^d} 
\\ & = \lim_r\frac{|\Omega_l \cap B(y,r)|}{\omega_{d}r^{d}} +\lim_r\frac{|\Omega_p \cap B(y,r)|}{\omega_{d}r^{d}} + \limsup_r\frac{|\Omega_i \cap B(y,r)|}{\omega_{d}r^{d}}
\\ & = 1 + \limsup_r\frac{|\Omega_i \cap B(y,r)|}{\omega_{d}r^{d}}
\end{split}
\end{equation}
which says that $y$ has density zero in $\Omega_i$.

The fact that (\ref{zerolimit}) holds is now a consequence of the geneneral fact

\begin{equation}
\limsup_{r\downarrow 0} \frac{\Hm( \Sigma_{ij} \cap B(x,r))}{\omega_{d-1}r^{d-1}} = 0\nonumber
\end{equation}
for $\Hm$-a.e. $x \in \Sigma_{ij}^c$.
\end{proof}

\begin{proof}[Proof of Lemma \ref{propertiesPOU}]
The argument for (1) can be found in \cite{Laux2016} in the case of two phases and without localization, i.e.\ with $\eta = 1$ and $N = 2$. For the sake of completeness, we provide the proof in our case. Upon splitting into the negative and positive part, we may assume $\eta \ge 0$. Clearly the only nonzero terms in the sum are those for wich $B_{m}^r \cap \Sigma_{ij} \neq \emptyset$. Fix such a ball: by definition there exists $y \in r\mathbf{Z}^d$ such that $B_{m}^r = B(y, 2r\sqrt{d})$. If $x \in \Sigma_{ij} \cap B_{m}^r$ then we have that $B(x, 2r\sqrt{d}) \subset B(y, 4r\sqrt{d})$, and by definition of $\mathcal{E}^r$ this yields

\begin{equation}
\Hm(B(x, 2r\sqrt{d})\cap \Sigma_{ij}) \le \frac{\omega_{d-1}}{2^d}(4r)^{d-1}\sqrt{d}^{d-1} = \frac{\omega_{d-1}}{2}(2r)^{d-1}\sqrt{d}^{d-1}.\nonumber
\end{equation}
Thus $x$ belongs to the set of points in $\Sigma_{ij} \cap B_{m}^r$ such that 

\begin{equation}\label{contradictionBound}
\frac{\Hm(B(x, 2r\sqrt{d})\cap \Sigma_{ij})}{\omega_{d-1}(2r\sqrt{d})^{d-1}} \le \frac{1}{2}.
\end{equation}

By De Giorgi's structure theorem the approximate tangent plane exists at every point $x \in \Sigma_{ij}$, thus (\ref{contradictionBound}) cannot hold when $r$ is small enough: moreover every point $x \in \Sigma_{ij}$ is contained in at most $c(2,d)$ balls, this means that

\begin{equation}\label{functionDomination}
\sum_{m} \mathbf{1}_{\left\{z \in B_{m}^r \cap \Sigma_{ij}:\ \frac{\Hm(B(x, 2r\sqrt{d})\cap \Sigma_{ij})}{\omega_{d-1}(2r\sqrt{d})^{d-1}} \le \frac{1}{2}\right\}}(x)\eta(x) \le c(2,d) \eta(x)
\end{equation}
and that the left hand side of (\ref{functionDomination}) converges to zero pointwise. By the dominated convergence theorem we get our claim.

Proof of (2). Upon splitting into the negative and positive part, we may assume $\eta \ge 0$. Given a point $x \in \Sigma_{lp}$, if $y \in r\mathbf{Z}^d$ is such that $x \in B(y, 2r\sqrt{d})$, then $B(y, 4r\sqrt{d}) \subset B(x, 6r\sqrt{d})$. Thus for any $1 \le i < j \le N$ with $(i,j) \neq (l,p)$ we have

\begin{align*}
\Hm(B(y,4r\sqrt{d}) \cap \Sigma_{ij}) & \le \Hm(B(x,6r\sqrt{d}) \cap \Sigma_{ij})
\\ &\le \frac{\omega_{d-1}}{2^d}(4r)^{d-1}\sqrt{d}^{d-1}
\end{align*}
provided $r$ is small enough, this follows from Lemma \ref{densInt}. Since $\mathcal{F}_2^r$ covers $\mathbf{R}^d$ we obtain that

\begin{equation}
x \in \bigcup_{m \in \mathbf{N}} B_{m}^r\nonumber
\end{equation}
for all $r$ small enough. In other words 

\begin{equation}
\lim_{r \downarrow 0}\sum_m \rho_{m}(x) \eta(x) = \eta(x)\nonumber
\end{equation}
pointwise on $\Sigma_{lp}$, and the argument of the limit on the right hand side is dominated by $\eta$. Thus we may once again appeal to the dominated convergence theorem and conclude the proof.

\end{proof}

\subsection{Consistency and Monotonicity}

The following results are essetially contained in \cite{Esedoglu2015} and \cite{Laux2016}, indeed the proofs may be adapted because we are assuming that $a_{ij}$ and $b_{ij}$ satisfy the triangle inequality.

\begin{lemma}\label{consistency}

Let $\chi \in L^1((0,T), BV([0,1)^d)^N)$ such that $\chi(\cdot, t) \in \mathcal{A}$ for a.e.\ $t$. Then

\begin{equation}
\lim_{h\downarrow 0} \int_0^T E_h(\chi)dt = \int_0^T E(\chi)dt.\nonumber
\end{equation}

Even more is true: for any $g \in C^{\infty}([0,1)^d)$ and any pair $1 \le i \neq j \le N$ we have

\begin{align}
\begin{aligned}
&\lim_{h\downarrow 0} \frac{1}{\sqrt{h}}\int_0^T\int g(x)(\chi_i(x,t)K_{ij}^h*\chi_j(x,t) + \chi_j(x,t)K_{ij}^h*\chi_j(x,t))dxdt 
\\ &= \int g(x)K_{ij}(z)|\nu_{ij}\cdot z|dzdxdt.
\end{aligned}\nonumber
\end{align}

\end{lemma}

\begin{lemma}\label{monotonicity}

For any $0 < h \le h_0$ we have

\begin{equation}
E_h(u) \ge \left( \frac{\sqrt{h_0}}{\sqrt{h} + \sqrt{h_0}}\right)^{d+1}E_{h_0}(u).\nonumber
\end{equation}

\end{lemma}

\subsection{Improved convergence of the energies}
The following Lemma is an improvement of the convergence of the energies, the proof of this result is contained, with minor modifications, in the paper \cite{Laux2016}, Corollary 3.7.

\begin{lemma}\label{noHtoH}

Let $u^h$ be a sequence of $[0,1]$-valued functions such that $u^h \to \chi$ in $L^1([0,1)^d\times (0,T))$ and 

\begin{equation}
\lim_{h\downarrow 0} \int_0^T E_h(u^h(t))dt = \int_0^T E(\chi(t))dt.
\end{equation}
Then we have that

\begin{align}
&\begin{aligned}
\lim_{h\downarrow 0} \frac{1}{\sqrt{h}}\int G_{\gamma}^h(z)|f^{\gamma}_h(z) - f^{\gamma}(z)|dz = 0,
\end{aligned}
\\ &\begin{aligned}
\lim_{h\downarrow 0} \frac{1}{\sqrt{h}}\int G_{\beta}^h(z)|f^{\beta}_h(z) - f^{\beta}(z)|dz = 0.\nonumber
\end{aligned}
\end{align}

Where we set

\begin{align}
&\begin{aligned}
f_h^{\gamma}(z) = \sum_{i,j}a_{ij}\int u_i^h(x,t)u_j^h(x-z,t)dxdt,\quad f^{\gamma}(z) = \sum_{i,j}a_{ij}\int \chi_i(x,t)\chi_j(x-z,t)dxdt,
\end{aligned}\nonumber
\\ &\begin{aligned}
f_h^{\beta}(z) = \sum_{i,j}b_{ij}\int u_i^h(x,t)u_j^h(x-z,t)dxdt,\quad f^{\beta}(z) = \sum_{i,j}b_{ij}\int \chi_i(x,t)\chi_j(x-z,t)dxdt.
\end{aligned}\nonumber
\end{align}
\end{lemma}

\subsection{Some inequalities}

Here we gather some elementary inequalities which are used frequently.

\begin{lemma} 

Let $a,b,a',b' \in \{0,1\}$, then the following inequalities hold.

\begin{align}
&\begin{aligned}\label{ineqAB}
|a - b| = a(1-b) + b(1-a)
\end{aligned}
\\ &\begin{aligned}\label{ineqLongAB}
|(a-a')_{+} - (b-b')_+| &+ |(a-a')_{-} - (b-b')_-|
\\ &\le |a-b| + |a'-b'|
\end{aligned}
\end{align}
\end{lemma}

\begin{proof}
The first identity follows by expanding $|a-b| = |a-b|^2$. The second one is proved in \cite{Laux2019}. For the sake of completeness, we reproduce the proof here. There are two cases. In the first one we have $(a-a')(b-b') \ge 0$ and we may assume upon replacing $(a,a',b,b')$ with $(-a,-a',-b,-b')$ that $(a-a')$ and $(b-b')$ are non-negative. Then (\ref{ineqLongAB}) reduces to

\begin{equation}
|(a-a') - (b-b')| \le |a-b| + |a'-b'|.\nonumber
\end{equation}
The second case is given by $(a-a')(b-b') \le 0$. By an argument as before we may assume $(a-a') \ge 0 \ge (b-b')$, thus (\ref{ineqLongAB}) reduces to

\begin{equation}
(a-a')+(b-b') \le |a-b| + |a'-b'|.\nonumber
\end{equation}
\end{proof}

\begin{lemma}
There exists a constant $C>0$ depending only on $N, \mathbb{A}, \mathbb{B}$ such that for any $v \in \mathcal{M}$

\begin{equation}\label{closeL2K}
\int|v - K^{h_0}*v|dx \le C\sqrt{h_0}E_h(v)\ \text{for all}\ h_0 \ge h.
\end{equation}

\end{lemma}

\begin{proof}
The proof of (\ref{closeL2K}) is aontained in the proof of Lemma 3 in \cite{Laux2019} for the two phases case when $K^h$ is the scaled version of the Gaussian with variance $1$. The same proof may be adapted to our setting because we still have monotonicity of the energy (Lemma \ref{monotonicity}) and we can prove essentially by the use of Jensen's inequality that

\begin{equation}
\int |v - K^h*v|dx \le C\sqrt{h}E_h(v).\nonumber
\end{equation}
\end{proof}
\section*{Acknowledgements}
This project has received funding from the Deutsche Forschungsgemeinschaft (DFG, German Research Foundation) under Germany's Excellence Strategy -- EXC-2047/1 -- 390685813.

\nocite{*}
\bibliography{biblio}{}
\bibliographystyle{plain}

\end{document}